\documentclass[11pt]{article}
\usepackage{amssymb,amsmath,latexsym}
\allowdisplaybreaks

\begin{document}

\begin{doublespace}

\def\1{{\bf 1}}
\def\ind{{\bf 1}}
\def\nn{\nonumber}

\def\sA {{\cal A}} \def\sB {{\cal B}} \def\sC {{\cal C}}
\def\sD {{\cal D}} \def\sE {{\cal E}} \def\sF {{\cal F}}
\def\sG {{\cal G}} \def\sH {{\cal H}} \def\sI {{\cal I}}
\def\sJ {{\cal J}} \def\sK {{\cal K}} \def\sL {{\cal L}}
\def\sM {{\cal M}} \def\sN {{\cal N}} \def\sO {{\cal O}}
\def\sP {{\cal P}} \def\sQ {{\cal Q}} \def\sR {{\cal R}}
\def\sS {{\cal S}} \def\sT {{\cal T}} \def\sU {{\cal U}}
\def\sV {{\cal V}} \def\sW {{\cal W}} \def\sX {{\cal X}}
\def\sY {{\cal Y}} \def\sZ {{\cal Z}}

\def\bA {{\mathbb A}} \def\bB {{\mathbb B}} \def\bC {{\mathbb C}}
\def\bD {{\mathbb D}} \def\bE {{\mathbb E}} \def\bF {{\mathbb F}}
\def\bG {{\mathbb G}} \def\bH {{\mathbb H}} \def\bI {{\mathbb I}}
\def\bJ {{\mathbb J}} \def\bK {{\mathbb K}} \def\bL {{\mathbb L}}
\def\bM {{\mathbb M}} \def\bN {{\mathbb N}} \def\bO {{\mathbb O}}
\def\bP {{\mathbb P}} \def\bQ {{\mathbb Q}} \def\bR {{\mathbb R}}
\def\bS {{\mathbb S}} \def\bT {{\mathbb T}} \def\bU {{\mathbb U}}
\def\bV {{\mathbb V}} \def\bW {{\mathbb W}} \def\bX {{\mathbb X}}
\def\bY {{\mathbb Y}} \def\bZ {{\mathbb Z}}
\def\R {{\mathbb R}} \def\RR {{\mathbb R}} \def\H {{\mathbb H}}
\def\n{{\bf n}} \def\Z {{\mathbb Z}}

\newcommand{\expr}[1]{\left( #1 \right)}
\newcommand{\cl}[1]{\overline{#1}}
\newtheorem{thm}{Theorem}[section]
\newtheorem{lemma}[thm]{Lemma}
\newtheorem{defn}[thm]{Definition}
\newtheorem{prop}[thm]{Proposition}
\newtheorem{corollary}[thm]{Corollary}
\newtheorem{remark}[thm]{Remark}
\newtheorem{example}[thm]{Example}
\numberwithin{equation}{section}
\def\ee{\varepsilon}
\def\qed{{\hfill $\Box$ \bigskip}}
\def\NN{{\mathcal N}}
\def\AA{{\mathcal A}}
\def\MM{{\mathcal M}}
\def\BB{{\mathcal B}}
\def\CC{{\mathcal C}}
\def\LL{{\mathcal L}}
\def\DD{{\mathcal D}}
\def\FF{{\mathcal F}}
\def\EE{{\mathcal E}}
\def\QQ{{\mathcal Q}}
\def\SS{{\mathcal S}}
\def\RR{{\mathbb R}}
\def\R{{\mathbb R}}
\def\L{{\bf L}}
\def\K{{\bf K}}
\def\S{{\bf S}}
\def\A{{\bf A}}
\def\E{{\mathbb E}}
\def\F{{\bf F}}
\def\P{{\mathbb P}}
\def\N{{\mathbb N}}
\def\eps{\varepsilon}
\def\wh{\widehat}
\def\wt{\widetilde}
\def\pf{\noindent{\bf Proof.} }
\def\pff{\noindent{\bf Proof} }
\def\cp{\mathrm{Cap}}

\title{\Large \bf Minimal thinness with respect to subordinate killed Brownian motions}

\author{{\bf Panki Kim}\thanks{This work was supported by the National Research Foundation of
Korea (NRF) grant funded by the Korea government (MEST) (NRF-2013R1A2A2A01004822)
}
\quad {\bf Renming Song\thanks{Research supported in part by a grant from
the Simons Foundation (208236)}} \quad and
\quad {\bf Zoran Vondra\v{c}ek}
\thanks{Research supported in part by the Croatian Science Foundation under the project 3526}
}

\date{}

\maketitle

\begin{abstract}
Minimal thinness is a notion that describes the smallness of a set
at a boundary point. In this paper, we provide  tests for
minimal thinness for a large class of subordinate killed
Brownian motions in
 bounded $C^{1,1}$ domains,
$C^{1,1}$ domains with compact complements and
domains above graphs of bounded $C^{1,1}$ functions.
\end{abstract}



\section{Introduction}\label{s:intro}

Let $X=(X_t,\P_x)$ be a Hunt process in an open set $D\subset \R^d$, $d\ge 2$.
Let $\partial_M D$ and $\partial_m D$ be the Martin and minimal Martin boundary
of $D$ with respect to $X$ respectively.
For any $z\in\partial_M D$,
we denote by $M^D(x,z)$  the Martin kernel of $D$ at $z$ with respect to $X$.
The family of all excessive functions for $X$ will be denoted by $\SS$. For a function
$v:D\to [0,\infty]$
 and a set $E\subset D$, the reduced function of $v$ on $E$ is
defined by $R^E_v=\inf\{s\in \SS:\, s\ge v \textrm{ on }E\}$ and its lower
semi-continuous regularization is denoted by $\wh{R}^E_v$.  A set $E\subset D$ is
said to be \emph{minimally thin} in $D$ at $z\in \partial_m D$ with respect to $X$
if $\wh{R}^E_{M^D(\cdot, z)}\neq M^D(\cdot, z)$, cf.~\cite{Fol}.
A probabilistic interpretation of minimal thinness is given in
terms of the process $X$ conditioned to die at $z\in \partial_m D$:
For any $z\in \partial_m D$, let $X^z=(X^z_t, \P^z_x)$ denote
the $M^D(\cdot, z)$-process, Doob's $h$-transform
of $X$ with $h(\cdot)=M^D(\cdot, z)$. The lifetime of $X^z$
will be denoted by $\zeta$. It is known (see \cite{KW}) that
$\lim_{t\uparrow\zeta}X^z_t=z$,   $\P^z_x$-a.\/s. For $E\subset D$,
let $T_E:=\inf\{t>0: X^z_t\in E\}$. It is proved in \cite[Satz 2.6]{Fol}
that a set $E\subset D$ is minimally thin at $z\in \partial_m D$ with respect to $X$
if and only if there exists $x\in D$ such that $\P^z_x(T_E<\zeta)\neq 1$.
This shows that minimal thinness is a concept describing smallness
of a set at a boundary point.

The history of minimal thinness
goes back to Lelong-Ferrand \cite{LF} who introduced this concept
in case of the half-space in the setting of classical potential theory.
Minimal thinness for general open sets
was developed in Na\"{i}m \cite{Nai}, while probabilistic
interpretation (in terms of Brownian motion) was given by
Doob (see e.g.~\cite{Doo}). Various versions of Wiener-type
criteria for minimal thinness were developed over the years
culminating in the work of Aikawa \cite{A} who, by using
the powerful concept of quasi-additivity of capacity,
established a criterion for minimal thinness for subsets of
NTA domains. For a good exposition of these results and
methods cf.~\cite[Part II, 7]{AE}. In case of a
$C^{1,1}$ domain $D\subset \R^d$, the finite part of the minimal Martin boundary
$\partial_m D$ coincides with the Euclidean boundary $\partial D$,
and Aikawa's criterion reads as follows: Let $E$ be a Borel subset
of $D$. If $E$ is minimally thin at $z\in \partial D$, then
\begin{equation}\label{e:aikawa-criterion}
\int_{E\cap B(z, 1)}|x-z|^{-d}\, dx <\infty\, .
\end{equation}
Conversely, if $E$ is the union of a subfamily of Whitney
cubes of $D$ and \eqref{e:aikawa-criterion} holds, then $E$
is minimally thin in $D$ at $z$.

Note that all works listed above pertain to
the classical potential theory related to Brownian motion. For more general
Hunt processes, although the general theory of minimal thinness was
developed by F\"ollmer already in 1969, see \cite{Fol}, until recently
no concrete criteria for minimal thinness were known. The first paper
addressing this question was \cite{KSV6} which dealt with
minimal thinness of subsets of the half-space for a large class of
subordinate Brownian motions. Quite general results
for a large class of symmetric L\'evy processes in $\kappa$-fat open
sets were obtained in \cite{KSV11}. The special case of a
$C^{1,1}$ open set $D$ was given in \cite[Corollary 1.5]{KSV11}.
We present here a slightly simplified version of the main result of \cite{KSV11}. Assume that $X$
is an isotropic L\'evy process in $\R^d$, $d\ge 2$, with
characteristic exponent $\Psi(x)=\Psi(|x|)$ satisfying the
following weak scaling condition: There exist constants
$0<\delta_1 \le \delta_2<1$ and $a_1,a_2>0$ such that
\begin{equation}\label{e:wsc}
a_1 \lambda^{2\delta_1}\Psi(t)\le \Psi(\lambda t)\le a_2
\lambda^{2\delta_2}\Psi(t)\, ,\qquad \lambda\ge 1, t\ge 1\, .
\end{equation}
We note that many subordinate Brownian motions, particularly all isotropic stable processes, satisfy the above
condition. Let $X^D$ be the process $X$ killed upon exiting
a $C^{1,1}$ open set $D$.  If a Borel set $E\subset D$ is
minimally thin in $D$ at $z\in \partial D$ with respect
to $X^D$, then \eqref{e:aikawa-criterion} holds true.
The converse is also true provided $E$ is the union of
a subfamily of Whitney cubes of $D$. Thus one obtains the
same Aikawa-type criterion for minimal thinness regardless
of the particular isotropic L\'evy process $X$ as long as $X$
satisfies the weak scaling condition \eqref{e:wsc}. This is a
somewhat surprising result. An explanation for this hinges on
sharp two-sided estimates for the Green function of $X^D$ which
imply that the singularity of the Martin kernel $M^D(x,z)$
near $z\in \partial D$ is of the order $|x-z|^{-d}$ for all such processes.

The purpose of this paper is to exhibit a large class of
(non-L\'evy) Markov processes
for which the Aikawa-type criterion for minimal thinness
depends on the particular process and is different from
\eqref{e:aikawa-criterion}. This class consists of subordinate
killed Brownian motions via subordinators having Laplace exponents
satisfying a certain weak scaling condition. Let us now precisely
formulate the setting and results.

Let $W=(W_t, \P_x)$ be a Brownian motion in $\R^d$, $d\ge 2$,
with transition density
$$
p(t,x,y)=(4\pi t)^{-\frac{d}{2}}\exp\left(-\frac{|x-y|^2}{4t}\right)\,
,\qquad t>0,\  x,y\in \R^d\, .
$$
Let $S=(S_t)_{t\ge 0}$ be an independent subordinator with Laplace
exponent $\phi:(0,\infty)\to (0,\infty)$, i.e.,
$\E[e^{-\lambda S_t}]=e^{-t\phi(\lambda)}$, $t\ge 0$, $\lambda>0$.
The process $X=(X_t, \P_x)$ defined by $X_t=W_{S_t}$, $t\ge 0$,
is called a subordinate Brownian motion. It is an isotropic L\'evy
process with characteristic exponent $\Psi(x)=\phi(|x|^2)$.
Let $D$ be an open subset of $\R^d$, and let $X^D$ be the process
$X$ killed upon exiting $D$. This process is known as a killed
subordinate Brownian motion. By reversing the order of subordination
and killing one obtains a different process. Assume from now on that
$D$ is a  domain (i.e.,~connected open set) in $\R^d$, and
let $W^D=(W^D_t, \P_x)$ be the Brownian motion $W$ killed upon
exiting $D$. The process $Y^D=(Y^D_t,\P_x)$ defined by
$Y^D_t=W^D_{S_t}$, $t\ge 0$, is called a subordinate killed
Brownian motion.
It is a Hunt process and its infinitesimal generator is given by
$-\phi(-\Delta_{|D})$ where $\Delta_{|D}$ is the Dirichlet Laplacian.

Recall that the Laplace exponent of a subordinator is a Bernstein function,
i.e.,~it has the representation
$$
\phi(\lambda)=b\lambda +\int_{(0,\infty)}(1-e^{-\lambda x})\, \mu(dx)\, ,
$$
with $b\ge 0$ and $\mu$ a measure on $(0,\infty)$ satisfying
$\int_{(0,\infty)}(1\wedge x)\, \mu(dx)<\infty$, which is called the L\'evy measure of $S$.
The potential measure of the subordinator $S$ is defined by
$
 U(A)=\int_0^\infty\P(S_t\in A)\,dt.
$
A Bernstein function $\phi$ is called a complete Bernstein function
if its L\'evy measure has a completely monotone density.
A Bernstein function $\phi$ is called a special Bernstein function if the
function $\lambda\mapsto \lambda/\phi(\lambda)$ is also a Bernstein function.
The function $\lambda\mapsto \lambda/\phi(\lambda)$ is called the conjugate
Bernstein function of $\phi$. It is well known that any complete Bernstein function
is a special Bernstein function. For this and other properties of
complete and special Bernstein functions, see \cite{SSV}.

In this the paper we will impose following assumptions:
\begin{itemize}
	\item[{\bf (A1)}]
	the potential measure of $S$ has a decreasing density $u$;
	\item[{\bf(A2)}]
	the L\' evy measure of $S$ is infinite and has a decreasing density $\mu$;
	\item[{\bf(A3)}] there exist constants $\sigma>0$, $\lambda_0>0$ and
$\delta \in (0, 1]$
 such that
\begin{equation*}
  \frac{\phi'(\lambda t)}{\phi'(\lambda)}\leq\sigma\, t^{-\delta}\ \text{ for all }\ t\geq 1\ \text{ and }\
 \lambda\geq\lambda_0\, .
\end{equation*}
\end{itemize}
Depending on whether our domain $D$ is bounded or unbounded, we will consider  the following two sets of conditions.

\begin{itemize}
	\item[{\bf(A4)}]
If $D$ is bounded and $d = 2$, we assume that there are $
\sigma_0>0$ and
$\delta_0 \in (0,2)$
such that
\begin{equation*}\label{e:new23}
  \frac{\phi'(\lambda t)}{\phi'(\lambda)}\geq
\sigma_0\,t^{-
\delta_0}\ \text{ for all
}\ t\geq 1\ \text{ and }\ \lambda\geq\lambda_0.
\end{equation*}

\item[{\bf(A5)}] If $D$ is bounded and $d = 2$, we assume that
\begin{equation*}
\int_{0}^{1}\frac{d \lambda}{\phi(\lambda)}<\infty.
\end{equation*}
\end{itemize}

\begin{itemize}
	\item[{\bf(A6)}]
If $D$ is unbounded then we assume that  $d  \ge 3$ and that there are
$\beta, \sigma_1 >0$ such that
\begin{align}\label{e:mas}
\frac{u(\lambda t)}{u(\lambda)} \ge  \sigma_1 t^{-\beta }
\quad \text{ for all } t\ge 1  \text{ and } \lambda>0\, .
\end{align}
\end{itemize}
Assumptions {\bf (A1)}--{\bf (A5)} were introduced and used in \cite{KM} and \cite{KM2}. It is easy to check that
if $\phi$ is a complete Bernstein function satisfying condition {\bf (H1)}: there exist $a_1, a_2>0$ and
$\delta_1, \delta_2\in (0,1)$ satisfying
$$
a_1 \lambda^{\delta_1}\phi(t)\le \phi(\lambda t)\le a_2 \lambda^{\delta_2}\phi(t)\, ,\qquad \lambda \ge 1, t\ge 1\, ,
$$
then {\bf (A1)}--{\bf (A4)} are automatically satisfied.
One of the reasons for adopting the more general setup above is to cover the case
of geometric stable and iterated geometric stable subordinators.
Suppose that $\alpha\in (0, 2)$ for $d \ge 2$ and that $\alpha\in (0, 2]$ for $d \ge 3$.
A geometric $(\alpha/2)$-stable subordinator is a subordinator
with Laplace exponent $\phi(\lambda)=\log(1+\lambda^{\alpha/2})$.
Let $\phi_1(\lambda):=\log(1+\lambda^{\alpha/2})$, and for $n\ge 2$,
$\phi_n(\lambda):=\phi_1(\phi_{n-1}(\lambda))$. A subordinator with
Laplace exponent $\phi_n$ is called an iterated geometric subordinator.
It is easy to check that the functions $\phi$ and $\phi_n$ satisfy
{\bf (A1)}--{\bf (A6)},
but they do not satisfy {\bf (H1)}.

Assumption {\bf (A1)} implies that $\phi$ is a special Bernstein
function, see, for instance, \cite[Theorem 5.1]{SV}.
Moreover,  {\bf (A3)} implies $b=0$,   {\bf (A2)} implies that $\mu((0, \infty))=\infty$, and
{\bf (A5)} is equivalent to the transience of  $X$.
In case $d\ge3$, $X$ is always transient.

Condition {\bf (A6)} is only assumed when $D$ is unbounded and can be restated as
\begin{equation}\label{e:mas2}
\frac{u(R)}{u(r)}\ge \sigma_1 \left(\frac{R}{r}\right)^{-\beta}\, ,\qquad 0<r\le R<\infty\, .
\end{equation}
Under {\bf (A1)}--{\bf (A3)}, the inequality in \eqref{e:mas2} is valid with $\beta=2-\delta$ whenever $0<r\le R\le 1$, (see \eqref{e:upper-estimate-u} and \eqref{e:lower-estimate-u} below).
So {\bf (A6)} is mainly a condition about the behavior of $u$ near
infinity. It follows easily from \cite{KSV8} that if $\phi$ is a complete Bernstein function satisfying,
in addition to {\bf (H1)}, also condition {\bf (H2)}: there exist $a_3, a_4>0$ and
$\delta_3, \delta_4\in (0,1)$ satisfying
$$
a_3 \lambda^{\delta_3}\phi(t)\le \phi(\lambda t)\le a_4 \lambda^{\delta_4}\phi(t)\, ,
\qquad \lambda \le 1, t\le 1\, ,
$$
then {\bf (A6)} is satisfied, see \cite[Corollary 2.4]{KSV8}.
There are plenty of examples of complete Bernstein functions which
satisfy {\bf (A6)} but not {\bf (H2)}. For any $m>0$ and $\alpha\in (0, 2)$, the function
$\phi(\lambda):=(\lambda+m^{2/\alpha})^{\alpha/2}-m$, the Laplace exponent of a
relativistic stable subordinator,
is such an example.

Recall that an open set $D$ in $\bR^d$ is said to be a (uniform) $C^{1,1}$
open set  if there exist a localization radius $R>0$ and
a constant $\Lambda>0$ such that for every $z\in\partial D$, there exist
a $C^{1,1}$-function $\psi=\psi_z: \bR^{d-1}\to \bR$ satisfying
$\psi (0)= 0$, $\nabla \psi (0)=(0, \dots, 0)$, $\| \nabla \psi
\|_\infty \leq \Lambda$, $| \nabla \psi (x)-\nabla \psi (w)|
\leq \Lambda |x-w|$, and an orthonormal coordinate system
$CS_z$ with its origin at $z$ such that
$$
B(z, R)\cap D=\{ y= (\wt y, \, y_d) \mbox{ in } CS_z:
|y|< R, y_d > \psi (\wt y) \}.
$$
The pair $(R, \Lambda)$ is called the characteristics
of the $C^{1,1}$ open set $D$.

Recall that an open set $D$ is said to satisfy the interior and exterior balls conditions with radius
$R_1$ if for every $z\in \partial D$, there exist $x\in D$ and $y\in \overline{D}^c$ such that
${\rm dist}(x, \partial D)=R_1$, ${\rm dist}(y, \partial D)=R_1$,
$B(x, R_1)\subset D$ and $B(y, R_1)\subset  \overline{D}^c$.
It is known, see \cite[Definition 2.1 and Lemma 2.2]{AKSZ}, that an open set
$D$ is a $C^{1,1}$ open set if and only if it satisfies the interior and exterior ball conditions.
By taking $R$ smaller if necessary, we will always assume a  $C^{1,1}$ open set
with characteristics $(R, \Lambda)$ also satisfies the interior and exterior balls conditions with the same radius $R$.

We can now state the main result of
this paper. By $\delta(x)$ we denote the distance of the point
$x\in D$ to the boundary $\partial D$.

\begin{thm}\label{t:main}
Assume that $\phi$ is a Bernstein function satisfying
{\bf (A1)}--{\bf (A6)}.
Let $D\subset \R^d$ be either a bounded $C^{1,1}$ domain, or
a $C^{1,1}$ domain with compact complement or
a domain above the graph of a bounded $C^{1,1}$ function.

\noindent
(1)
If $E$ is  minimally thin in $D$  at $z\in\partial D$ with
respect to $Y^D$,
then
\begin{equation}\label{e:main}
\int_{E\cap B(z, 1)} \frac{\delta(x)^2\phi(\delta(x)^{-2})
\phi'(|x-z|^{-2})}{|x-z|^{d+4}\phi(|x-z|^{-2})^2}\, dx < \infty\, .
\end{equation}

\noindent
(2)
Conversely, if $E$ is the union of a subfamily of Whitney
cubes of $D$ and \eqref{e:main} holds true, then $E$ is
minimally thin in $D$  at $z\in\partial D$ with respect to $Y^D$.
\end{thm}

Since minimal thinness is defined for points in the minimal Martin boundary, the first step in proving this theorem is the identification of the finite part of the (minimal) Martin boundary of $D$
with its Euclidean boundary. In case of a bounded Lipschitz domain,
special subordinator $S$, and $d\ge 3$, this was accomplished in \cite[Theorem 4.3]{SV06}
(see also \cite[Theorem 5.84]{SV}). The method employed in
\cite{SV06, SV} heavily depended on the fact that the semigroup of the killed Brownian
motion $W^D$ in a bounded Lipschitz domain $D$ is intrinsically ultracontractive
which implies that all excessive functions with respect to $W^D$ are purely excessive.
In fact, \cite{SV06} proves that there is 1-1 correspondence between the cone of excessive
(respectively non-negative harmonic) functions of $W^D$ and the cone of excessive
(respectively non-negative harmonic) functions of $Y^D$, thus
allowing an easy transfer of many results valid for $W^D$ to results for $Y^D$.
In case of an unbounded domain, the semigroup of $W^D$ is no longer intrinsically ultracontractive and the method from \cite{SV06} cannot be used to identify the finite part of the (minimal) Martin boundary of $D$ with its Euclidean boundary.

In the case of killed subordinate Brownian motions, one of the main
tools used in identifying the (minimal) Martin boundary of a (possibly)
unbounded open set is the boundary Harnack principle.

In the
present case of subordinate killed Brownian motions, the
boundary Harnack principle is not yet available. As a substitute for the
boundary Harnack principle, we first establish sharp two-sided estimates
on the Green functions of subordinate killed Brownian motions in any $C^{1,1}$ domain
with compact complement or any domain above the graph of a bounded $C^{1,1}$ function.
This is done in Section \ref{s:gfmke}, see Theorems \ref{t:gfe} and \ref{t:gfe:nn}.
In Section \ref{s:mke}, by using some ideas from \cite{SV06},  we then show that the Martin kernel $M^D_Y(\cdot, \cdot)$
can be extended from $D\times D$ to $D\times \overline{D}$, cf.~Proposition \ref{p:extension-M}.
By using sharp two-sided estimates of the Green function, we subsequently establish in
Theorems \ref{t:mke} and \ref{t:mke2} sharp two-sided estimates
for the Martin kernel $M^D_Y(x,z)$, $x\in D$, $z\in \partial D$.
The remaining part of the section is devoted to proving that the finite part of the (minimal)
Martin boundary of $D$ can be identified with its Euclidean boundary in case
$D$ is either a bounded $C^{1,1}$ domain,
a $C^{1,1}$ domain with compact complement or a domain above the graph of a bounded $C^{1,1}$ function.
We note that in case of a bounded $C^{1,1}$ domain
(and under the assumptions {\bf (A1)}--{\bf (A5)})
this gives an alternative proof of some of the results form \cite{SV06}. Results of Sections \ref{s:gfmke} and \ref{s:mke} might be of independent interest.

Having identified the finite part of the (minimal) Martin boundary with the Euclidean boundary,
we can follow the method developed by Aikawa, cf.~\cite{A} and \cite[Part II, 7]{AE}, which
was also used in \cite{KSV11}, to prove Theorem \ref{t:main}.
One of the main ingredients of this method is the quasi-additivity of
the capacity related to the process $Y^D$, see Proposition
\ref{p:quasi-additivity}. This depends on the construction of a measure
comparable to the capacity which relies on an appropriate Hardy's inequality.
The first result on minimal thinness is a criterion given in Proposition
\ref{p:minthin-criterion-1} stating that a subset $E$ of $D$ is minimally
thin at $z\in \partial D$ (with respect to $Y^D$) if and only if
$\sum_{n=1}^{\infty}R^{E_n}_{M^D_Y(\cdot, z)}(x_0)<\infty$; here
$E_n=E\cap \{x\in D:\, 2^{-n-1}\le |x-z|<2^{-n}\}$ and $x_0\in D$
a fixed point. The proof of this general result depends
on an inequality relating the Green function and the Martin kernel of $Y^D$,
cf. Corollary \ref{c:U-M}. The inequality itself hinges on
sharp two-sided estimates of the Green function of $Y^D$ (cf. Theorems \ref{t:gfe} and \ref{t:gfe:nn})
and sharp two-sided
estimates of the Martin kernel (cf. Theorems \ref{t:mke} and \ref{t:mke2}).
With the quasi-additivity of capacity and the criterion for
minimal thinness from Proposition \ref{p:minthin-criterion-1} in hand,
it is rather straightforward to complete the proof of Theorem \ref{t:main}.

As an application of Theorem \ref{t:main}, we derive an analogue to a
criterion in the classical setting for minimal thinness in the half-space $\H$ of
a set below the graph of a Lipschitz function
$f:\R^{d-1}\to [0,\infty)$. In the classical case and
the case of killed subordinate Brownian motions in the half-space studied in \cite{KSV11}, the
criterion states that the set $A=\{(\wt{x}, x_d)\in
\H:\, 0<x_d\le f(\wt{x})\}$ is minimally
thin at 0 if and only
if $\int_{\{|\wt{x}|<1\}}f(\wt{x})|\wt{x}|^{-d}\, d\wt{x} <\infty$.
For the subordinate killed Brownian motion $Y^D$ the criterion
depends on the underlying Bernstein function $\phi$ and says
that $A$ is minimally thin at 0 if and only if
$$
\int_{\{|\wt{x}|<1\}}\frac{f(\wt{x})^3\phi(
f(\wt{x})^{-2})\phi'(|\wt{x}|^{-2})}{|\wt{x}|^{d+4}\phi(|\wt{x}|^{-2})^2}\, d\wt{x} <\infty\, ,
$$
see Proposition \ref{p:dahlberg2} and Remark \ref{r:dahlberg2-2} for the precise statement.

Finally, we give some examples. We first look at
three processes related to the stable process:
(1) $X^D$ -- the isotropic $\alpha$-stable process killed upon
exiting $D$, (2) $Y^D$ -- the subordinate killed Brownian motion
in $D$ with $(\alpha/2)$-stable subordinator, and
(3) $Z^D$ -- the censored $\alpha$-stable process in $D$.
Following \cite{MV} we briefly indicate how to prove criteria
for minimal thinness for the censored process, and then
compare minimal thinness of a given set with respect to
these processes and the index of stability $\alpha$.
Roughly, minimal thinness for $Z^D$ implies minimal
thinness for $X^D$ which in turn implies minimal
thinness for $Y^D$, see Corollary \ref{c:comparison}
for the precise statement. We also show that the converse does not hold.
At the end of Section \ref{s:stable}, we give some examples
related to subordinate killed Brownian motions
via geometric stable subordinators.

Organization of the paper: In the next section we give some
preliminaries on Bernstein functions satisfying conditions {\bf (A1)}--{\bf (A5)} and on the subordinate killed Brownian motion $Y^D$
and its relation to the killed subordinate Brownian motion.
In Section \ref{s:gfmke} we prove sharp two-sided estimates
for the Green function and the jumping kernel of $Y^D$.
In Section \ref{s:mke}
we identify the finite part of the (minimal) Martin boundary with
the Euclidean boundary
and  give
sharp  two-sided estimates
on the Martin kernel of $Y^D$.
We continue in Section \ref{s:quasi} with the proof of the
quasi-additivity of the capacity.
Results about minimal thinness
are proved in Section \ref{s:mtf}. The paper concludes with
criteria for minimal thinness with respect to processes related to the stable case,
and with respect to subordinate killed Brownian motions via
geometric stable subordinators.

In this paper, we use the letter $c$, with
or without subscripts, to denote a constant, whose value may change from one
appearance to another. The notation $c(\cdot, \dots, \cdot)$ specifies the dependence of
the constant. The dependence of the constants on the domain $D$
(including the dimension $d$) and
the Bernstein function $\phi$ will not be explicitly mentioned.
For any two positive functions $f$ and $g$,
$f\asymp g$ means that there is a positive constant $c\geq 1$
so that $c^{-1}\, g \leq f \leq c\, g$ on their common domain of
definition.
We will use ``$:=$" to denote a
definition, which is read as ``is defined to be".
For $a, b\in \bR$,
$a\wedge b:=\min \{a, b\}$ and $a\vee b:=\max\{a, b\}$.


\section{Preliminaries}\label{s:preliminaries}

In this section we first collect several properties of Bernstein
functions and then collect some results on the subordinate killed Brownian motion $Y^D$ and its relation to the killed subordinate Brownian motion $X^D$.

\begin{lemma}\label{l:properties-of-bf}
\begin{itemize}
\item[(a)] For every Bernstein function $\phi$,
\begin{equation}\label{e:wlsc-substitute}
1 \wedge  \lambda\le \frac{\phi(\lambda t)}{\phi(t)} \le 1 \vee \lambda\, ,\quad \textrm{for all }t>0, \lambda>0\, .
\end{equation}
\item[(b)]
If $\phi$ is a special Bernstein function, then
$\lambda \mapsto \lambda^2\phi'(\lambda)$ and $\lambda\mapsto
\lambda^2\frac{\phi'(\lambda)}{\phi(\lambda)^2}$ are increasing functions.
Furthermore, for any $\gamma>2$, $\lim_{\lambda\to 0} \lambda^\gamma\frac{\phi'(\lambda)}{\phi(\lambda)^2}=0$.
\item[(c)]
If $\phi$ is a special Bernstein function, then for
every $d\ge 2$, $\gamma\ge 2$, $\lambda>0$, $b\in (0,1]$ and $a\in [1,\infty)$ it holds that
\begin{equation}\label{e:phi-prime}
\frac{b}{a^{d+\gamma+1} \lambda^{d+\gamma}}\frac{\phi'(\lambda^{-2})}{\phi(\lambda^{-2})^2}\le \frac1{t^{d+\gamma}}\frac{\phi'(t^{-2})}{\phi(t^{-2})^2} \le \frac{a}{b^{d+\gamma+1} \lambda^{d+\gamma}}\frac{\phi'(\lambda^{-2})}{\phi(\lambda^{-2})^2}\, ,\quad\textrm{for all }t\in [b\lambda, a\lambda]\, .
\end{equation}
\end{itemize}
\end{lemma}

Part (a) is well known, part (b) is proved in \cite[Lemma 4.1]{KM}, and part (c) can be proved in the same way as \cite[Corollary 2.2]{KM2} where the proof is given for $\gamma=2$. We will frequently use all three properties of the lemma, often without explicitly mentioning it.

Let $W$ be a Brownian motion in $\R^d$, $D\subset \R^d$ a domain, and $W^D$ a Brownian motion killed upon exiting $D$. We denote by $p^D(t,x,y)$, $t>0$, $x,y\in D$, the transition densities of $W^D$, and by $(P_t^D)_{t\ge 0}$ the corresponding semigroup.
Let $S$ be a subordinator independent of the Brownian motion $W$.
Let $Y^D_t=W^D_{S_t}$ be the
corresponding subordinate killed Brownian motion in $D$.
The process $Y^D$ is a symmetric Hunt process, cf. \cite{SV08}. We will use
$(\EE^D, \DD(\EE^D))$ to denote the Dirichlet form associated with $Y^D$.
The killing measure of $\EE^D$
has a density $\kappa_D$ given by the formula
\begin{equation}\label{e:killing-function-Y}
\kappa_D(x)=\int_{(0,\infty)}(1-P^D_t 1(x))\, \mu(dt)\, ,\quad x\in D\, .
\end{equation}
It follows from the general theory of Dirichlet forms that for every $v\in \DD(\EE^D)$ it holds that
\begin{equation}\label{e:dirichlet-form-Y}
\EE^D(v,v)\ge \int_D v(x)^2 \kappa_D(x)\, dx\, .
\end{equation}

Let $(R^D_t)_{t\ge 0}$ be the transition semigroup of $Y^D$.
We will need to compare this semigroup with the semigroup of the killed
subordinate Brownian motion. Recall that $X_t=W_{S_t}$ is the
subordinate Brownian motion and $(X^D_t)_{t\ge 0}$
is the subprocess of $X$ killed upon exiting $D$.
Let $(Q^D_t)_{t\ge 0}$ denote the transition semigroup of $X^D$.
It is well known, cf.~\cite[Proposition 3.1]{SV08}, that
$(R^D_t)_{t\ge 0}$ is subordinate to $(Q^D_t)_{t\ge 0}$ in
the sense that
\begin{align}
\label{e:subord}
R^D_t f(x)\le Q^D_t f(x)  \quad \text{ for all Borel }
f:D\to [0,\infty) \text{ all } t\ge 0 \text{ and all } x\in D.
\end{align}

Let $j_X(x)$ denote the density of the L\'evy measure of the process $X$. Then
$$
j_X(x)=\int_{(0,\infty)}p(t,x,0)\, \mu(dt)=\int_{(0,\infty)}(4\pi t)^{-\frac{d}{2}}
\exp\left(-\frac{|x|^2}{4t}\right)\, \mu(dt)\, .
$$
Clearly, $j_X$ is a continuous function of $x$ on $\R^d \setminus \{0\}$ and radial
(that is, $j_X(x)=j_X(|x|)$).
Let $\kappa_D^X$ denote the killing function of $X^D$. Then
\begin{equation}\label{e:kappa-X}
\kappa_D^X(x)=\int_{D^c}j_X(x-y)\, dy\, , \quad   x\in D\, ,
\end{equation}
and $\kappa_D^X$ is a continuous function of $x\in D$.

\begin{lemma}\label{l:transition-densities}
For any open set $D \subset \R^d$,
\begin{equation}\label{e:killing-functions}
 \kappa_D^X(x)\le \kappa_D(x)\, ,\quad
\text{for almost all }x\in D\, .
\end{equation}
\end{lemma}
\pf
Using \eqref{e:subord},
the Lemma  follows from the argument of \cite[Proposition 3.2]{SV03}. \qed

Assume $\phi$ is a Bernstein function satisfying
{\bf (A1)} so that the potential measure of $S$ has a decreasing density $u(t)$. Then the Green function of the subordinate killed Brownian motion $Y^D$, denoted by $U^D(x,y)$, $x,y \in D$,  is given by the formula
\begin{equation}\label{e:green-function-skbm}
U^D(x,y)=\int_0^{\infty}p^D(t,x,y)u(t)\, dt=\int_0^{\infty}r^D(t,x,y)\, dt\, , \quad x,y \in D\, .
\end{equation}
Similarly, the Green function of $X$, denoted by $G_X(x,y)$, $x,y\in \R^d$, is given by
\begin{equation}\label{e:green-function-sbm}
G_X(x,y)=\int_0^{\infty}p(t,x,y)u(t)\, dt\, , \quad x,y \in \R^d\, .
\end{equation}
Since $p^D(t,x,y)\le p(t,x,y)$ for all $x,y\in D$, we see from \eqref{e:green-function-skbm} and \eqref{e:green-function-sbm} that
\begin{equation}\label{e:U-less-GX}
U^D(x,y)\le G_X(x,y) \, ,\quad \textrm{for all }x,y\in D\, .
\end{equation}

Assume now that $\phi$ is a Bernstein function satisfying
{\bf (A1)}--{\bf (A5)} and let $S$ be a subordinator with Laplace exponent $\phi$.
The potential density  $u(t)$ of $S$ satisfies the following two estimates:
\begin{align}
u(t)\le (1-2e^{-1})^{-1}\frac{\phi'(t^{-1})}{t^2\phi(t^{-1})^2}\, ,\quad t>0\, ,\label{e:upper-estimate-u}
\end{align}
and, for every $M>0$ there exists $c_1=c_1(M)>0$ such that
\begin{align}
u(t)\ge c_{1} \frac{\phi'(t^{-1})}{t^2\phi(t^{-1})^2}\, ,\quad 0<t\le
M\, .\label{e:lower-estimate-u}
\end{align}
For the upper estimate see \cite[Lemma A.1]{KM},
and for the lower \cite[Proposition 3.4]{KM}

The density $\mu(t)$ of the L\'evy measure of $S$ satisfies the following two estimates:
\begin{align}
\mu(t)\le (1-2e^{-1})^{-1}t^{-2}\phi'(t^{-1})\, ,\quad t>0\, ,
\label{e:upper-estimate-mu}
\end{align}
and, for every $M>0$ there exists $c_2=c_2(M)>0$ such that
\begin{align}
\mu(t)\ge c_2 t^{-2}\phi'(t^{-1})\, ,\quad 0<t\le
M\, .\label{e:lower-estimate-mu}
\end{align}
For the upper estimate see \cite[Lemma A.1]{KM},
and for the lower \cite[Proposition 3.3]{KM}.

Recall that $G_X(x,y)$ denotes the Green function of the subordinate
Brownian motion $X_t=W_{S_t}$. When $d \ge 3$ we have  that there exists $c_3>0$ such that
\begin{equation}\label{e:free-gfe}
G_X(x,y)\le c_3 \frac{\phi'(|x-y|^{-2})}
{|x-y|^{d+2}\phi(|x-y|^{-2})^2}\, ,\quad x,y\in \R^d\, .
\end{equation}
This can be proved by following the proof of \cite[Lemma 3.2(b)]{KSV8} using \eqref{e:upper-estimate-u}
and \cite[Lemma 4.1]{KM}.
Moreover, by
\cite[Proposition 4.5]{KM} we have the following two-sided inequality:
For every $d \ge 2$ and $M>0$,
there exists $c_4=c_4(M)>1$ such that
\begin{equation}\label{e:free-gfe-reverse}
c_4^{-1}
\frac{\phi'(|x-y|^{-2})}
{|x-y|^{d+2}\phi(|x-y|^{-2})^2} \le G_X(x,y)\le c_4
\frac{\phi'(|x-y|^{-2})}
{|x-y|^{d+2}\phi(|x-y|^{-2})^2}\, ,
\quad |x-y|\le M\, .
\end{equation}

The L\'evy density of $X$ also has the following two-sided estimates by \cite[Proposition 4.2]{KM}:
For every $M>0$ there exists $c_5=c_5(M)>0$ such that
\begin{equation}\label{e:levy-density-estimate}
c_5^{-1}r^{-d-2}\phi'(r^{-2})\le j_X(r)\le c_5 r^{-d-2}\phi'(r^{-2})\, ,\quad r\in (0, M]\, .
\end{equation}
Thus, by using Lemma \ref{l:properties-of-bf}(a) and (c), for every $M>0$,
\begin{equation}\label{e:doubling-condition}
j_X(r)\le c j_X(2r)\, ,\quad r\in (0, M]\, .
\end{equation}

\section{Kernel estimates on subordinate killed Brownian motion}\label{s:gfmke}

In this section we assume that
$D\subset \R^d$ is either a bounded $C^{1,1}$ domain, or
a $C^{1,1}$ domain with compact complement or a
domain above the graph of a bounded  $C^{1,1}$
function. We assume that the $C^{1, 1}$ characteristics
of $D$ is $(R, \Lambda)$.

Recall that $(P^D_t)_{t\ge 0}$ denotes the transition semigroup of the killed Brownian motion $W^D$ and $p^D(t,x,y)$, $t>0$, $x,y\in D$, is the corresponding transition density.
It is known that $p^D(t,x,y)$ satisfies the following short-time
estimates (cf.~\cite{Zh, Zh2, Son}):
For any $T>0$, there exist positive constants
$c_1, c_2, c_3, c_4$
such that for any $t\in (0, T]$ and any $x,y\in D$,
\begin{equation}\label{e:upper-bound-for-p2}
p^D(t,x,y)\le c_1\left(\frac{\delta_D(x)}{\sqrt{t}} \wedge 1 \right)
\left(\frac{\delta_D(y)}{\sqrt{t}} \wedge 1 \right)
\, t^{-d/2}\exp\left(-\frac{c_2 |x-y|^2}{t}\right)\, ,
\end{equation}
\begin{equation}\label{e:lower-bound-for-p}
p^D(t,x,y)\ge c_3\left(\frac{\delta_D(x)}{\sqrt{t}} \wedge 1 \right)
\left(\frac{\delta_D(y)}{\sqrt{t}} \wedge 1 \right)
\, t^{-d/2}\exp\left(-\frac{c_4 |x-y|^2}{t}\right)\, .
\end{equation}
Thus, by the semigroup property and \eqref{e:upper-bound-for-p2}, we get
there exist positive constants $c_5, c_6, c_7, c_8$ such that
for every $t>3$
\begin{align*}
&p^D(t,x,y)
=\int_D \int_D p^D(1,x,z) p^D(t-2,z,w)p^D(1,w,y) dzdw\\
&\le c_5\left(\delta(x)\wedge 1 \right)\left(\delta(y)\wedge 1 \right)\\
&\quad \times
\int_D \int_D
\,\exp\left(-c_6|x-z|^2\right)
\,  (t-2)^{-d/2}\exp\left(-\frac{c_6|z-w|^2}{t-2}\right)
\,\exp\left(-{c_6 |w-y|^2}\right) dzdw\\
&\le c_5\left(\delta(x)\wedge 1 \right)\left(\delta(y)\wedge 1 \right)\\
&\quad \times
\int_{\R^d} \int_{\R^d}
\,\exp\left(-c_6|x-z|^2\right)
\,  (t-2)^{-d/2}\exp\left(-\frac{c_6|z-w|^2}{t-2}\right)
\,\exp\left(-{c_6 |w-y|^2}\right) dzdw\\
&\le c_7
\left(\delta(x)\wedge 1 \right)\left(\delta(y)\wedge 1 \right)
t^{-d/2}\exp\left(-\frac{c_8|x-y|^2}{t}\right).
\end{align*}
Combining this with \eqref{e:upper-bound-for-p2}, we have that
there exist positive constant $c_9, c_{10}$
such that for all $t>0$ and any $x,y\in D$,
\begin{equation}\label{e:upper-bound-for-p}
p^D(t,x,y)\le c_9 \left(\frac{\delta(x)}{\sqrt{t} \wedge 1} \wedge 1 \right)
\left(\frac{\delta(y)}{\sqrt{t}\wedge 1} \wedge 1 \right)
\, t^{-d/2}\exp\left(-\frac{ c_{10} |x-y|^2}{t}\right)\, .
\end{equation}
We will use the following bound several times: By the change of variables
$s=c|x-y|^2/t$, for every $c>0$ and $a  \in \R$, we have
\begin{eqnarray}\label{e:UD-lower-2}
&&\int_0^{|x-y|^2}
\left(\frac{\delta(x)}{\sqrt{t}} \wedge 1 \right)
\left(\frac{\delta(y)}{\sqrt{t}} \wedge 1 \right)
t^{-a/2}\exp\left(
-\frac{c |x-y|^2}{t}\right) \, dt \nonumber\\
&=& \int_{ c}^{\infty}
\left(\frac{\sqrt{s/c}\, \delta(x)}{|x-y|} \wedge 1 \right)
\left(\frac{\sqrt{s/c}\, \delta(y)}{|x-y|} \wedge 1 \right)
\left(\frac{c|x-y|^2}{s}\right)^{-a/2}e^{-s}
\frac{c|x-y|^2}{s^2}\, ds\nonumber\\
&\ge& c^{1-(a/2)} \left(\frac{\delta(x)}{|x-y|} \wedge 1 \right)
\left(\frac{\delta(y)}{|x-y|} \wedge 1 \right)
|x-y|^{-a+2} \int_{ c}^{\infty} s^{a/2-2}
e^{-s}\, ds.
\end{eqnarray}

Our first goal is to obtain sharp two-sided estimates on $U^D$.
Under stronger assumptions on the Laplace exponent $\phi$ such
estimates were given in \cite[Theorem 5.91]{SV} for bounded $D$.
In the remainder of this section
$\phi$ is a Bernstein function satisfying {\bf (A1)}--{\bf (A5)}.
We first consider the case $|x-y| \le M$.
\begin{thm}\label{t:gfe}
For every
$M>0$,
there exists a constant $c=c(M)\ge 1$ such that
for all $x,y \in D$ with $|x-y|\le M$,
\begin{eqnarray}
&&c^{-1}\left(\frac{\delta(x)\delta(y)}{|x-y|^2}\wedge 1
\right)\frac{\phi'(|x-y|^{-2})}
{|x-y|^{d+2}\phi(|x-y|^{-2})^2}\nonumber\\
&& \le U^D(x,y)
\le c \left(\frac{\delta(x)\delta(y)}{|x-y|^2}\wedge
1\right)\frac{\phi'(|x-y|^{-2})}
{|x-y|^{d+2}\phi(|x-y|^{-2})^2}\, .\label{e:gfe}
\end{eqnarray}
\end{thm}

\pf {\it Upper bound}:
It follows from \eqref{e:U-less-GX} and \eqref{e:free-gfe-reverse} that there exists a constant $c_1>0$ such that for
all $x,y \in D$ with $ |x-y|\le M$,
\begin{equation}\label{e:upper-UD-1}
U^D(x,y)\le G_X(x,y)\le c_1 \frac{\phi'(|x-y|^{-2})}
{|x-y|^{d+2}\phi(|x-y|^{-2})^2}\, .
\end{equation}

Let $c_2$ be the constant $c_{10}$ in \eqref{e:upper-bound-for-p}.
Since $t \to \frac{\phi'(t^{-1})}{\phi(t^{-1})^2}$ is
 increasing,  using \eqref{e:upper-estimate-u} we have that for $r>0$,
\begin{align}\label{e:uewq11}
 I_1(r)&:=\int_0^{r^2}t^{-d/2-1}\exp{\left(-\tfrac{c_2r^2}{t}\right)}
u(t)\,dt\leq c_3\int_0^{r^2}
t^{-d/2-1}\exp{\left(-\tfrac{c_2r^2}{t}\right)}
t^{-2}\frac{\phi'(t^{-1})}{\phi(t^{-1})^2}\,dt\nonumber\\&\leq
c_3\frac{\phi'(r^{-2})}{\phi(r^{-2})^2}\int_0^{r^2}t^{-\frac{d}{2}-
3}\exp{\left(-\tfrac{c_2r^2}{t}\right)}\,dt
=c_4r^{-d-4}\frac{\phi'(r^{-2})}{\phi(r^{-2})^2}\int_{{c_2}}^\infty
t^{\frac{d}{2}+1}e^{-t}\,dt\,.
\end{align}
On the other hand, since $u$ is decreasing, using \eqref{e:upper-estimate-u}
we have that  for $r>0$,
\begin{align}\label{e:uewq21}
 I_2(r)&:=\int_{r^2}^\infty
t^{-d/2-1}u(t)\,dt    \leq  u(r^2)
\int_{r^2}^\infty t^{-d/2-1} dt\nonumber\\
&\leq c_5 r^{-4} \frac{\phi'(r^{-2})}{\phi(r^{-2})^2} \int_{r^2}^\infty t^{-d/2-1} dt\le c_6
r^{-d-4}\frac{\phi'(r^{-2})}{\phi(r^{-2})^2}\, .
\end{align}
It follows from \cite[Lemma 4.4]{KM} that
\begin{align}\label{e:uewq3}
 L:=\int_{(2M)^2}^\infty
t^{-d/2}u(t)\,dt   <\infty.
\end{align}
Thus from \eqref{e:green-function-skbm}, \eqref{e:upper-bound-for-p} and  \eqref{e:uewq11}--\eqref{e:uewq3},
we have that, for $|x-y|\le M$,
\begin{align*}
&U^D(x,y)\,=\,\int_0^{\infty}p^D(t,x,y)u(t)\, dt\\
\le&
\int_0^{|x-y|^2}p^D(t,x,y)u(t)\, dt +\int_{|x-y|^2}^{(2M)^2}p^D(t,x,y)u(t)\, dt+\int_{(2M)^2}^\infty p^D(t,x,y)u(t)\, dt\\
\le& c_7
\int_0^{|x-y|^2}t^{-d/2-1}\delta(x)\delta(y)\exp
\left(-\frac{c_2|x-y|^2}{t}\right)u(t)\, dt \\
&
+c_7
\int_{|x-y|^2}^{(2M)^2} t^{-d/2-1}\delta(x)\delta(y)
u(t)\, dt
+c_7
\int_{(2M)^2}^{\infty} t^{-d/2}\delta(x)\delta(y)u(t)\, dt
\nonumber \\
\le &c_7 \delta(x)\delta(y) \big(I_1(|x-y|)+ I_2(|x-y|)+L\big) \le c_8\frac{\delta(x)\delta(y)}{|x-y|^2}\,
\frac{\phi'(|x-y|^{-2})}
{|x-y|^{d+2}\phi(|x-y|^{-2})^2}\, .
\end{align*}
In the last inequality above we use the fact that
$r\to r^{-d-4}\frac{\phi'(r^{-2})}{\phi(r^{-2})^2}$ is a decreasing function and is thus bounded from below
by a positive constant on $(0, M^2]$.
Together with \eqref{e:upper-UD-1} this gives the upper bound
in \eqref{e:gfe}.

\noindent
{\it Lower bound}:
Since $u$ is decreasing and $|x-y| \le M$, by \eqref{e:lower-bound-for-p} and
\eqref{e:lower-estimate-u},
\begin{eqnarray}\label{e:UD-lower-1}
&&U^D(x,y)\ge c_{9}\int_0^{|x-y|^2}
\left(\frac{\delta(x)}{\sqrt{t}} \wedge 1 \right)
\left(\frac{\delta(y)}{\sqrt{t}} \wedge 1 \right)
 t^{-d/2}\exp\left(
-\frac{c_{10} |x-y|^2}{t}\right) u(t)\, dt \nonumber\\
&&\ge c_{9}u(|x-y|^2)\int_0^{|x-y|^2}
\left(\frac{\delta(x)}{\sqrt{t}} \wedge 1 \right)
\left(\frac{\delta(y)}{\sqrt{t}} \wedge 1 \right)
t^{-d/2}\exp\left(
-\frac{c_{10} |x-y|^2}{t}\right) \, dt\nonumber\\
&&\ge c_{11}\frac{\phi'(|x-y|^{-2})}
{|x-y|^{4}\phi(|x-y|^{-2})^2} \int_0^{|x-y|^2}
\left(\frac{\delta(x)}{\sqrt{t}} \wedge 1 \right)
\left(\frac{\delta(y)}{\sqrt{t}} \wedge 1 \right)
t^{-d/2}\exp\left(
-\frac{c_{10} |x-y|^2}{t}\right) \, dt.\nonumber
\end{eqnarray}
By combining this with \eqref{e:UD-lower-2}
we arrive at
\begin{align*}
U^D(x,y)&\ge c_{12}
\left(\frac{\delta(x)}{|x-y|} \wedge 1 \right)
\left(\frac{\delta(y)}{|x-y|} \wedge 1 \right)
\frac{\phi'(|x-y|^{-2})}
{|x-y|^{d+2}\phi(|x-y|^{-2})^2}\nonumber\\
&\asymp
\left(\frac{\delta(x)\delta(y)}{|x-y|^2}\wedge 1\right)
\frac{\phi'(|x-y|^{-2})}
{|x-y|^{d+2}\phi(|x-y|^{-2})^2}\, .
\end{align*}
\qed

We now assume $d \ge 3$ and consider our two types of unbounded $C^{1,1}$ domains and give different estimates for $U^D$.

If $D\subset \R^d$ is a
domain above the  graph of a bounded $C^{1,1}$
function, then it follows from \cite{Zh, Son}
that there exist positive constants $c_1, c_2, c_3$ and $c_4$
such that for any $t\in (0, \infty)$ and any $x,y\in D$,
\begin{equation}\label{e:upper-bound-for-pp2}
p^D(t,x,y)\le c_1\left(\frac{\delta(x)}{\sqrt{t}} \wedge 1 \right)
\left(\frac{\delta(y)}{\sqrt{t}} \wedge 1 \right)
\, t^{-d/2}\exp\left(-\frac{c_2 |x-y|^2}{t}\right)\, ,
\end{equation}
\begin{equation}\label{e:lower-bound-for-pp2}
p^D(t,x,y)\ge c_3\left(\frac{\delta(x)}{\sqrt{t}} \wedge 1 \right)
\left(\frac{\delta(y)}{\sqrt{t}} \wedge 1 \right)
\, t^{-d/2}\exp\left(-\frac{c_4 |x-y|^2}{t}\right)\, .
\end{equation}
Clearly for $a>2$,
\begin{align}
\label{e:compu2-1}
\int_{|x-y|^2}^\infty \left(\frac{\delta(x)}{\sqrt{t}} \wedge 1 \right)
\left(\frac{\delta(y)}{\sqrt{t}} \wedge 1 \right)
\, t^{-a/2} dt
\le  \frac{2}{a-2}
 \left(\frac{\delta(x)}{|x-y|} \wedge 1 \right)
\left(\frac{\delta(y)}{|x-y|} \wedge 1 \right)\frac{1}
{|x-y|^{a-2}}\, .
\end{align}
By the change of variables
$s=|x-y|^2/t$ and the inequality
$$
\left(\frac{\sqrt{s}\delta(x)}{|x-y|} \wedge 1 \right)
\le  \sqrt{s}\left(\frac{\delta(x)}{|x-y|} \wedge 1 \right),
\qquad s\ge 1,
  $$
it is easy to see  that for $a \in \R$ and $b>0$,
there exist a constant $c=c(a, b)>0$ such that
\begin{align}
\label{e:compu2-2}
\int_0^{|x-y|^2} \left(\frac{\delta(x)}{\sqrt{t}} \wedge 1 \right)
\left(\frac{\delta(y)}{\sqrt{t}} \wedge 1 \right)
\, t^{-a/2}\exp\left(-\frac{b |x-y|^2}{t}\right) dt \nonumber
\\
\le  c
 \left(\frac{\delta(x)}{|x-y|} \wedge 1 \right)
\left(\frac{\delta(y)}{|x-y|} \wedge 1 \right)\frac{1}
{|x-y|^{a-2}}\, .
\end{align}

If $D\subset \R^d$ is a
$C^{1,1}$ domain with compact complement, then it follows from \cite{Zh2} that
there exist positive constants $c_5, c_6, c_7$ and $c_8$
such that for any $t\in (0, \infty)$ and any $x,y\in D$,
\begin{equation}\label{e:upper-bound-for-pp3}
p^D(t,x,y)\le c_5\left(\frac{\delta(x)}{\sqrt{t} \wedge 1} \wedge 1 \right)
\left(\frac{\delta(y)}{\sqrt{t} \wedge 1} \wedge 1 \right)
\, t^{-d/2}\exp\left(-\frac{c_6 |x-y|^2}{t}\right)\, ,
\end{equation}
\begin{equation}\label{e:lower-bound-for-pp3}
p^D(t,x,y)\ge c_7\left(\frac{\delta(x)}{\sqrt{t} \wedge 1} \wedge 1 \right)
\left(\frac{\delta(y)}{\sqrt{t} \wedge 1} \wedge 1 \right)
\, t^{-d/2}\exp\left(-\frac{c_8 |x-y|^2}{t}\right)\, .
\end{equation}

Clearly
for $a>2$,
\begin{align}
&\int_{|x-y|^2}^\infty \left(\frac{\delta(x)}{\sqrt{t}\wedge 1} \wedge 1 \right)
\left(\frac{\delta(y)}{\sqrt{t}\wedge 1} \wedge 1 \right)
\, t^{-a/2} dt\nonumber\\
&  \le  \frac{2}{a-2}
 \left(\frac{\delta(x)}{|x-y|\wedge 1} \wedge 1 \right)
\left(\frac{\delta(y)}{|x-y|\wedge 1} \wedge 1 \right)\frac{1}
{|x-y|^{a-2}}\, . \label{e:compu3-1}
\end{align}

By the change of variables
$s=|x-y|^2/t$ and the inequalities
$$
\left(\frac{\delta(x)}{|x-y|\wedge 1} \wedge 1 \right) \le \left(\frac{\delta(x)}{(|x-y|/\sqrt{s})\wedge 1} \wedge 1 \right)
  \le
  \sqrt{s}\left(\frac{\delta(x)}{|x-y|\wedge 1} \wedge 1 \right),   \quad  s\ge 1 ,
  $$
it is easy to see that for $a \in \R$ and $b>0$, there exists a constant $c=c(a, b)>0$ such that
\begin{align}
\int_0^{|x-y|^2} \left(\frac{\delta(x)}{\sqrt{t}\wedge 1} \wedge 1 \right)
\left(\frac{\delta(y)}{\sqrt{t}\wedge 1} \wedge 1 \right)
\, t^{-a/2}\exp\left(-\frac{b |x-y|^2}{t}\right) dt \nonumber\\ \le  c
 \left(\frac{\delta(x)}{|x-y|\wedge 1} \wedge 1 \right)
\left(\frac{\delta(y)}{|x-y|\wedge 1} \wedge 1 \right)\frac{1}
{|x-y|^{a-2}} \label{e:compu3-2}
\end{align}
and
\begin{align}
\label{e:compl3}
\int_0^{|x-y|^2} \left(\frac{\delta(x)}{\sqrt{t}\wedge 1} \wedge 1 \right)
\left(\frac{\delta(y)}{\sqrt{t}\wedge 1} \wedge 1 \right)
\, t^{-a/2}\exp\left(-\frac{b |x-y|^2}{t}\right) dt \nonumber \\
\ge  c^{-1}
 \left(\frac{\delta(x)}{|x-y|\wedge 1} \wedge 1 \right)
\left(\frac{\delta(y)}{|x-y|\wedge 1} \wedge 1 \right)\frac{1}
{|x-y|^{a-2}}\, .
\end{align}

\begin{thm}\label{t:gfe:nn}
Suppose that $d \ge 3$ and that $\phi$ is a
Bernstein function satisfying {\bf (A1)}--{\bf (A3)} and {\bf (A6)}.
(1) Let $D\subset \R^d$ be a
domain above the graph of a bounded $C^{1,1}$
function.
There exists a constant $c_1\ge 1$ such that for all $x,y \in D$,
\begin{align*}
&c_1^{-1}\left(\frac{\delta(x)}{|x-y|} \wedge 1 \right)
\left(\frac{\delta(y)}{|x-y|} \wedge 1 \right)\frac{u(|x-y|^2)}
{|x-y|^{d-2}} \le U^D(x,y)\nn\\
&\qquad \le c_1 \left(\frac{\delta(x)}{|x-y|} \wedge 1 \right)
\left(\frac{\delta(y)}{|x-y|} \wedge 1 \right)\frac{u(|x-y|^2)}
{|x-y|^{d-2}}\, .
\end{align*}
(2) Let $D\subset \R^d$ be a
$C^{1,1}$ domain with compact complement. There exists a constant $c_1\ge 1$ such that
for all $x,y \in D$,
\begin{align*}
&c_1^{-1}\left(\frac{\delta(x)}{|x-y|\wedge 1} \wedge 1 \right)
\left(\frac{\delta(y)}{|x-y|\wedge 1} \wedge 1 \right)\frac{u(|x-y|^2)}
{|x-y|^{d-2}} \le U^D(x,y) \nn\\
&\qquad \le c_1 \left(\frac{\delta(x)}{|x-y|\wedge 1} \wedge 1 \right)
\left(\frac{\delta(y)}{|x-y|\wedge 1} \wedge 1 \right)\frac{u(|x-y|^2)}
{|x-y|^{d-2}}\, .
\end{align*}
\end{thm}

\pf
We give the proof of (2) first.

\noindent
{\it Upper bound}:
Using \eqref{e:mas} and the fact  $u$ is decreasing, we have
from \eqref{e:upper-bound-for-pp3} that
\begin{align*}
&U^D(x,y)\,=\,\int_0^{\infty}p^D(t,x,y)u(t)\, dt\\
\le& c_1
\int_0^{\infty}\left(\frac{\delta(x)}{\sqrt{t}\wedge 1} \wedge 1 \right)
\left(\frac{\delta(y)}{\sqrt{t}\wedge 1} \wedge 1 \right)t^{-d/2}\exp
\left(-\frac{c_2|x-y|^2}{t}\right)u(t)\, dt
\nonumber \\
\le& c_3 |x-y|^{2\beta} u(|x-y|^2)
\int_0^{|x-y|^2}\left(\frac{\delta(x)}{\sqrt{t}\wedge 1} \wedge 1 \right)
\left(\frac{\delta(y)}{\sqrt{t}\wedge 1} \wedge 1 \right)t^{-\beta-d/2}\exp
\left(-c_2\frac{|x-y|^2}{t}\right)\, dt \\
&
+c_1 u(|x-y|^2)
\int_{|x-y|^2}^{\infty} \left(\frac{\delta(x)}{\sqrt{t}\wedge 1} \wedge 1 \right)
\left(\frac{\delta(y)}{\sqrt{t}\wedge 1} \wedge 1 \right)t^{-d/2}\, dt.
\end{align*}
Together with \eqref{e:compu3-1}--\eqref{e:compu3-2} we obtain the upper bound.

\noindent
{\it Lower bound}:
Since $u$ is decreasing,  by
\eqref{e:lower-bound-for-pp3}
\begin{align}\label{e:UD-lower-111}
U^D(x,y)\ge &c_4\int_0^{|x-y|^2}
\left(\frac{\delta(x)}{\sqrt{t}\wedge 1} \wedge 1 \right)
\left(\frac{\delta(y)}{\sqrt{t}\wedge 1} \wedge 1 \right)
 t^{-d/2}\exp\left(
-\frac{c_5 |x-y|^2}{t}\right) u(t)\, dt \nonumber\\
\ge &c_4u(|x-y|^2)\int_0^{|x-y|^2}
\left(\frac{\delta(x)}{\sqrt{t}\wedge 1} \wedge 1 \right)
\left(\frac{\delta(y)}{\sqrt{t}\wedge 1} \wedge 1 \right)
t^{-d/2}\exp\left(
-\frac{c_5 |x-y|^2}{t}\right) \, dt\, .
\end{align}
Combining \eqref{e:UD-lower-111} and \eqref{e:compl3}
we arrive at
\begin{align*}
U^D(x,y)&\ge c_6
\left(\frac{\delta(x)}{|x-y|\wedge 1} \wedge 1 \right)
\left(\frac{\delta(y)}{|x-y|\wedge 1} \wedge 1 \right)
\frac{u(|x-y|^2)}
{|x-y|^{d-2}}\, . \end{align*}

Using \eqref{e:UD-lower-2} and \eqref{e:upper-bound-for-pp2}--\eqref{e:compu2-2},
instead of \eqref{e:upper-bound-for-pp3}--\eqref{e:UD-lower-111},
the proof of (1) is similar to (2).
\qed

\begin{prop}\label{p:U-continuous}
The Green function $U^D$ is  jointly continuous in the extended sense, hence jointly lower
semi-continuous, on $D\times D$.
\end{prop}
\pf Let $x,y\in D$, $x\neq y$, and set $\eta=|x-y|/2$. Let $(x_n,y_n)_{n\ge 1}$ be a sequence in $D\times D$ converging to $(x,y)$ and  assume that $|x_n-y_n|\ge \eta$. For every $t>0$, $\lim_{n\to \infty}p^D(t,x_n,y_n)=p^D(t,x,y)$. Moreover
$$
p^D(t,x_n,y_n)\le (4\pi t)^{-d/2}\exp\left(-\frac{|x_n-y_n|^2}{4t}\right)\le (4\pi t)^{-d/2}\exp\left(-\frac{\eta^2}{4t}\right)\, .
$$
Since the process $X$ is transient, we have that
$$
\int_0^{\infty}(4\pi t)^{-d/2}\exp\left(-\frac{\eta^2}{4t}\right)u(t)\, dt<\infty\, .
$$
Now it follows from the bounded convergence theorem that
$$
\lim_{n\to \infty}U^D(x_n,y_n)=\lim_{n\to \infty}\int_0^{\infty}p^D(t,x_n,y_n)u(t)\, dt=\int_0^{\infty}p^D(t,x,y)u(t)\, dt =U^D(x,y)\, .
$$
On the other hand,
from Theorem \ref{t:gfe} we get that
$$
\lim_{(x_n,y_n)\to (x,x)}U^D(x_n,y_n)=+\infty =U^D(x,x)\, .
$$
Thus
$U^D$ is jointly continuous in the extended sense, and therefore jointly lower semi-continuous. \qed

We now recall a result from analysis (see \cite[Theorem 1, p. 167]{Stein}):
Any open set $D\subset\R^d$ is the union of a family $\{Q_j\}_{j \in \N}$
of closed cubes, with sides all parallel to the axes, satisfying
the following properties:
(i)  $\mathrm{int}(Q_j) \cap\, \mathrm{int} (Q_k)=\emptyset$, $j\neq k$;
(ii) for any $j$, $\mathrm{diam}(Q_j)\le \mathrm{dist}(Q_j,
\partial D)\le 4  \mathrm{diam}(Q_j)$, where $\mathrm{dist}
(Q_j, \partial D)$ denotes the Euclidean distance between
$Q_j$ and $\partial D$. The family $\{Q_j\}_{j \in \N}$ above
is called a Whitney decomposition of $D$ and the $Q_j$'s are
called Whitney cubes (of $D$). We will use $x_j$ to denote
the center of the cube $Q_j$. For each cube $Q_j$ let $Q_j^*$
denote the interior of the double of $Q_j$.

\begin{corollary}\label{c:whitney-green}
(i) For every $M>0$ there exists a constant $c_1=c_1(M)\ge 1$ such that for all
Whitney cubes $Q_j$ whose diameter is less than $M$,
\begin{align}\label{e:whitney-green}
c_1^{-1} U^D(x',y)\le U^D(x,y) \le c_1 U^D(x',y)\, ,
\end{align}
for all $x,x'\in Q_j$  and all $y\in D \setminus Q_j^*$ with ${\rm dist}(y, Q_j)<M$.

\noindent
(ii) For every $M>0$ there exists a constant $c_2=c_2(M) >0 $ such that for all cubes
$Q_j$ whose diameter is less than $M$ and all $x,x'\in Q_j$, it holds that
\begin{align}\label{e:whitney-green-2}
U^D(x,x')\ge c_2 G_X(x,x')\, .
\end{align}
\end{corollary}
\pf (i) From the geometry of Whitney cubes it is easy to see
that there exists a constant $c\ge 1$ such that for every cube
$Q_j$ it holds that
\begin{eqnarray*}
& &c^{-1}\delta(x)\le \delta(x_j)\le c \delta(x)\, ,\quad
\text{for all }x\in Q_j\, ,\\
& &c^{-1}|x-y|\le |x_j-y|\le c|x-y|\, ,\quad \text{for all }
x\in Q_j \text{ and all }y\in D\setminus Q_j^*\, .
\end{eqnarray*}
Together with Theorem \ref{t:gfe} and Lemma \ref{l:properties-of-bf}(c),
these estimates imply that
$$
U^D(x,y)\asymp U^D(x_j,y)\, ,\quad \text{for all }x\in Q_j
\text{ and all }y\in D\setminus Q_j^*\text{ with dist}(y, Q_j)<M\, ,
$$
with a constant independent of $Q_j$. This clearly implies
the statement of the corollary.

\noindent (ii) If $x,x'\in Q_j$, then $|x-x'|\le \mathrm{diam }
(Q_j)\le \mathrm{dist}(Q_j,\partial D)\le \delta(x)\wedge \delta(x') \wedge (4M)$.
Thus it follows from \eqref{e:gfe} and \eqref{e:free-gfe-reverse} that
$$
U^D(x,x')\,\ge\,c_1\,
\frac{\phi'(|x-x'|^{-2})}
{|x-x'|^{d+2}\phi(|x-x'|^{-2})^2}\,\ge\,
c_2\,G_X(x,x')\, .
$$
\qed

Let $J^D(x,y)$ be the jumping density of $Y^D$ defined by
$$J^D(x,y) =\int_0^{\infty}p^D(t,x,y)\mu(t)\, dt.$$
Clearly $J^D(x,y) \le j_X(|x-y|)$, $x, y \in D$.

Using \eqref{e:upper-estimate-mu}, \eqref{e:lower-estimate-mu}, \eqref{e:levy-density-estimate} and
the fact that $t^2\phi'(t)$ is increasing (see Lemma \ref{l:properties-of-bf}(b)),
the proof of the next proposition is very similar to that of Theorem
\ref{t:gfe}.

\begin{prop}
\label{p:J^D(z,y)}
For every $M>0$,
there exists a constant $c=c(M)\ge 1$ such that such that for all $x,y \in D$ with $|x-y|\le M$,
$$
c^{-1}\left(\frac{\delta(x)\delta(y)}{|x-y|^2}\wedge 1
\right)\frac{\phi'(|x-y|^{-2})}
{|x-y|^{d+2}} \le J^D(x,y)
\le c \left(\frac{\delta(x)\delta(y)}{|x-y|^2}\wedge
1\right)\frac{\phi'(|x-y|^{-2})}
{|x-y|^{d+2}}\, .
$$
\end{prop}

For any open subset $B$ of $D$, let
$U^{D,B}(x,y)$ be the Green function of $Y^D$ killed upon exiting $B$.
We define the Poisson kernel
\begin{equation}\label{PK}
K^{D,B}(x,y)\,:=  \int_{B}
U^{D,B}(x,z) J^D(z,y) dz, \qquad (x,y) \in
B \times  (D \setminus \overline{B}).
\end{equation}
Using the L\'{e}vy system for $Y^D$, we know that for every
open subset $B$ of $D$ and  every $f \ge 0$ on $D \setminus \overline{B}$ and $x \in B$,
\begin{equation}\label{newls}
\E_x\left[f(Y^D_{\tau_B});\,Y^D_{\tau_B-} \not= Y^D_{\tau_B}  \right]
=\int_{D \setminus \overline{B}} K^{D,B}(x,y)f(y)dy.
\end{equation}

\begin{lemma}
\label{l:KSV7prop4.7}
For every $M>0$, there exists $c=c(M)>0$ such that for
any ball $B(x_0, r)\subset D$ of radius $r\in (0, 1]$, we have
for all $(x, y)\in B(x_0, r)\times (D \setminus  \overline{B(x_0, r)})$ with $|x-y|\le M$,
\begin{eqnarray}\label{e:KSV7-4.7}
K^{D, B(x_0, r)}(x, y)\,\le \,c \,\delta(y)\frac{\phi'((|y-x_0|-r)^{-2})}
{(|y-x_0|-r)^{d+3}}    \phi(r^{-2})^{-1}.
\end{eqnarray}
\end{lemma}
\pf
Let $B=B(x_0, r)$.
Since $U^{D,B}(x,y) \le G_X(x,y)$,
\eqref{PK} and Proposition \ref{p:J^D(z,y)} imply that for  every $(x,y) \in
B \times  (D \setminus \overline{B})$ with $|x-y|\le M$,
\begin{align}
\label{e:KBE}
K^{D,B}(x,y)
&\le  \int_{B}G_X(x,z) J^D(z,y) dz\nn \\
&\le c_1(M) \int_{B}
G_X(x,z)\left(\frac{\delta(z)\delta(y)}{|z-y|^2}\wedge
1\right)\frac{\phi'(|z-y|^{-2})}
{|z-y|^{d+2}} dz \nn\\
&\le c_1(M) \delta(y)\int_{B}
G_X(x,z)\frac{\phi'(|z-y|^{-2})}
{|z-y|^{d+3}} dz.
\end{align}
Since $|z-y| \ge |y-x_0|-r$ and $t \to {t^{-d-3}}\phi'(t^{-2})$ is decreasing (see Lemma \ref{l:properties-of-bf}(b)),
 \begin{align}\label{e:KBE1}
\int_{B} G_X(x,z)\frac{\phi'(|z-y|^{-2})}{|z-y|^{d+3}} dz
&\le \frac{\phi'((|y-x_0|-r)^{-2})}
{(|y-x_0|-r)^{d+3}} \int_{B} G_X(x,z)dz\nn\\
&\le \frac{\phi'((|y-x_0|-r)^{-2})}
{(|y-x_0|-r)^{d+3}}\int_{B(0, 2r)} G_X(0, z)dz.
\end{align}
By \eqref{e:free-gfe-reverse}, we have
\begin{align}
\int_{B(0, 2r)} G_X(0, z)dz &\le
c_2 \int_{B(0, 2r)} |z|^{-d-2}\frac{\phi'(|z|^{-2})}
{\phi(|z|^{-2})^2}  dz=c_2  \int_0^{2r} r^{-3}
\frac{\phi'(r^{-2})}
{\phi(r^{-2})^2}  dr\nn\\
& \le 2^{-1}c_3\phi(2^{-1} r^{-2})^{-1} \le  2c_4\phi(r^{-2})^{-1}.\label{e:KBE2}
\end{align}
Combining  \eqref{e:KBE}--\eqref{e:KBE2}, we have proved the proposition.
\qed


\section{Martin boundary and Martin kernel estimates}\label{s:mke}

In this section
we assume that
$D\subset \R^d$ is either a bounded $C^{1,1}$ domain, or
a $C^{1,1}$ domain with compact complement or
a domain above the graph of a bounded $C^{1,1}$
function. We assume that the $C^{1,1}$ characteristics
of $D$ is $(R, \Lambda)$.

Denote by $\widetilde{Y}^D$ the subordinate  killed Brownian motion via
a subordinator with Laplace exponent $\lambda/\phi(\lambda)$.
Let $\widetilde \mu(dt)$ be the L\'evy measure of the
(possibly killed) subordinator with
Laplace exponent $\lambda/\phi(\lambda)$, the conjugate Bernstein function of $\phi(\lambda)$.
Since
$\mu((0, \infty))=\infty$, we also have  $\widetilde \mu((0, \infty))=\infty$,
$$
\frac{\lambda}{\phi(\lambda)}=u(\infty)+\int_0^\infty (1-e^{-\lambda t}) \widetilde \mu(dt)
$$
and
\begin{align}
\label{e:umu}
u(t)=\widetilde
\mu((t, \infty))
 +u(\infty).
\end{align}
(See \cite[Corollary 5.5]{SV} and the paragraph after it.)

Denote by $( \widetilde R^D_t)_{t\ge 0}$ the transition semigroup of $\widetilde Y^D$ and
by $\widetilde U^D$  the potential operator of $\widetilde Y^D$.
For any function $f$ which is excessive for $W^D$
we define an  operator $\widetilde V^D$ by
$$
\widetilde V^D f(x)=u(\infty) f(x)+ \int_{(0,\infty)}(f(x)-P^D_t f(x))\,
 \widetilde \mu(dt)\, ,\quad x\in D\, .
$$
Let $G^D( x, y)=\int_0^\infty p^D(t,x,y)dt$ be the Green function of $W^D$.

\begin{lemma}\label{l:VGD0}
For any $x,y \in D$, we have
$$U^D(x,y) =\widetilde V^D (G^D( \cdot, y))(x).$$
\end{lemma}
\pf
By the semigroup property,  for every $s>0$,
\begin{align*}
&G^D( x, y)
=\int_0^\infty p^D(t,x,y)dt =\int_0^s p^D(t,x,y)dt+ \int_0^\infty p^D(t+s,x,y)dt\\
&=\int_0^s p^D(t,x,y)dt+ P^D_s\int_0^\infty p^D(t,\cdot,y) (x)dt=
\int_0^s p^D(t,x,y)dt+P_s^DG^D (\cdot,y)(x).
\end{align*}
Thus
\begin{align}
\label{e:uvg1}
\int_{(0,\infty)}(G^D( x, y)-P_s^DG^D (\cdot,y)(x))\widetilde \mu(ds)
=\int_{(0,\infty)}\int_0^s p^D(t,x,y)dt\widetilde \mu(ds).
\end{align}
Using \eqref{e:umu}
we see that
\begin{align*}
\widetilde V^D (G^D( \cdot, y))(x)&=
u(\infty) G^D( x, y) +
\int_{(0,\infty)}\int_0^s p^D(t,x,y)dt\widetilde \mu(ds)\\
&=u(\infty) G^D( x, y) +
\int_{0}^\infty\widetilde \mu((t, \infty)) p^D(t,x,y)dt \\
&=u(\infty) G^D( x, y) +
\int_{0}^\infty (u(t) -u(\infty) ) p^D(t,x,y)dt =U^D(x,y).
\end{align*}
\qed

Note that according to the pointwise version of the Bochner subordination formula one can regard $-\widetilde V$ as the generator of $\widetilde{Y}^D$. This provides an intuitive explanation of Lemma  \eqref{l:VGD0}, namely $V^D U^D (\cdot, y)=V^D \widetilde{V}^D G^D(\cdot, y)=-\Delta G^D (\cdot, y)=-\delta_y$.

Fix a point $x_0\in D$ and define the Martin kernel with respect to $Y^D$ based at $x_0$ by
\begin{align}\label{e:mk}
M^D_Y(x, y):=\frac{U^D(x, y)}{U^D(x_0, y)}, \qquad x, y\in D,~y\neq x_0.
\end{align}
We will establish some relation between the Martin kernel for $Y^D$ and the Martin kernel for $W^D$.
Define the Martin kernel with respect to $W^D$ based at $x_0$ by
\begin{align}\label{e:mkW}
M^D(x, y):=\frac{G^D(x, y)}{G^D(x_0, y)}, \qquad x, y\in D,~y\neq x_0.
\end{align}
Since $D$ is a $C^{1,1}$ domain, for each $z\in \partial D$ there exists the limit
$$
M^D(x,z):=\lim_{y\to z} M^D(x,y)\, .
$$

In the next lemma, we extend  \cite[Lemma 5.82]{SV} by including our two types of unbounded $C^{1,1}$ domains and the case $d=2$ for bounded $C^{1,1}$ domains.

\begin{lemma}\label{t5.1}
If
$(y_j)_{j\ge 1}$
is a sequence of points in $D$ such
that $\lim_{j \to \infty} y_j=z\in \partial D$,
then for each $t>0$ and
each $x\in D$,
$$
\lim_{j\to \infty} P^D_t \left(\frac{G^D(\cdot,y_j)}{G^D(x_0, y_j)}\right)(x)
= P^D_t(M^D(\cdot,z))(x)\, .
$$
\end{lemma}

\pf
Recall that the $C^{1, 1}$ characteristics of $D$ is $(R, \Lambda)$.
Fix $x \in D$ and let $R_1:=(R\wedge |x_0-z| \wedge |x-z|)/4$.
We assume all $y_j$ are in $B(z, R_1/2) \cap D$.
For any $r\in (0, R_1]$, there exists a ball $B(A_r(z), r/2) \subset  D\cap B(z, r)$.
It is well known  (see \cite[page 140]{A1} and \cite[Theorem 7.1]{JK}) that
there exist $c_1, \beta>0$ such that for any $r\in (0, R_1]$
and any $(y, w) \in D \cap B(z, r) \times (D\setminus B(z, 2r))$,
\begin{equation}\label{e:martin-estimate-boundary}
|M^D(w, y)-M^D(w, z)|\le c_1 M^D(w, A_r(z))\left(\frac{|y-z|}{r}\right)^\beta\, .
\end{equation}

 Let $g(w)=|w|^{-d+2}$ be the Newtonian kernel when
$d\ge 3$ and be the logarithmic kernel $g(x)=\left(\log\frac{1}{|x|}\right)\vee 1$ when $d=2$.
Using the estimate of $p^D(t,x,y)$ in \eqref{e:upper-bound-for-p2} and the Green function estimates of Brownian motion, we have the following estimates:
for every $t>0$ there exists a constant $c_2=c_2(t, \delta(x), R_1)>0$ such that
\begin{align}
\label{pDMWD1}
p^D(t,x,y) M^D(y, z)\le c_2 g(y-z) \quad \forall y \in B(z, R_1) \cap D,
\end{align}
\begin{align}
\label{pDMWD2}
p^D(t,x,y)M^D(y, y_j)\le c_2 g(y-y_j) \quad \forall y \in B(y_j, R_1) \cap D.
\end{align}
In fact,  since
$$\left(\frac{\delta(y)}{|y-y_j|}\wedge \frac{\delta(y)}{\delta(y_j)} \right) \le 2,$$
for $d \ge3$,
\begin{align*}
&p^D(t,x,y)M^D(y, y_j) \le c_3(t) \delta(x) \delta(y) \frac{G^D(y, y_j)}{\delta(y_j)}\\
&\le c_4(t, \delta(x)) \frac{\delta(y)}{\delta(y_j)} \left(\frac{\delta(y)}{|y-y_j|}\wedge 1 \right)
\left(\frac{\delta(y_j)}{|y-y_j|}\wedge 1 \right) |y-y_j|^{-d+2}\\
&\le c_4(t, \delta(x)) \frac{\delta(y)}{\delta(y_j)} \left(\frac{\delta(y_j)}{|y-y_j|}\wedge 1 \right) |y-y_j|^{-d+2}\\
&\le c_4(t, \delta(x))\left(\frac{\delta(y)}{|y-y_j|}\wedge \frac{\delta(y)}{\delta(y_j)} \right)  |y-y_j|^{-d+2}
\le 2 c_4(t, \delta(x))  |y-y_j|^{-d+2}.
\end{align*}
This proves \eqref{pDMWD2} for $d\ge 3$, and  by letting $y_j \to z$, we get \eqref{pDMWD1} for $d\ge 3$.
The proofs of  \eqref{pDMWD1} and \eqref{pDMWD2} for $d=2$ are similar.

The inequalities \eqref{pDMWD1} and \eqref{pDMWD2} imply that for every $r \le R_1$ and sufficiently large $j$,
\begin{equation}\label{e:bm3g2}
\int_{D \cap B(z, r)} p^D(t,x,y)(M^D(y, y_j)+M^D(y, z))dy \le
2c_2 \int_{B(0, 2r)}  g(y) dy.
\end{equation}
Given $\eps>0$, choose $0<r_1\le R_1$ small such that
$ \int_{B(0, 2r_1)}  g(y) dy< \eps/(4c_2)$.
For $y\in D \setminus  B(z, r_1)$, by
\eqref{e:martin-estimate-boundary} we get that
\begin{align}\label{e:sdafe1}
|M^D(y, y_j)-M^D(y, z)|&\le c_2M^D(y, A_{r_1}(z))\left(\frac{|y_j-z|}{r_1}\right)^\beta.
\end{align}
Therefore, using the fact that $y \to M^D(y, z)$ is excessive for $W^D$, for every large $j$
\begin{align*}
&|P^D_t \left(\frac{G^D(\cdot,y_j)}{G^D(x_0, y_j)}\right)(x)- P^D_t(M^D(\cdot,z))(x)|\\
\le&
\int_{D \cap B(z, r_1)} p^D(t,x,y)(M^D(y, y_j)+M^D(y, z))dy +
\int_{D \setminus  B(z, r_1)} p^D(t,x,y)|M^D(y, y_j)-M^D(y, z)|dy \\
\le & \eps/2+ c_2\left(\frac{|y_j-z|}{r_1}\right)^\beta \int_{D}
P_t^DM^D(\cdot, A_{r_1}(z))(y)dy \le
 \eps/2+ c_2\left(\frac{|y_j-z|}{r_1}\right)^\beta M^D(x, A_{r_1}(z)) \le
\eps.
\end{align*}
\qed

Using the previous lemma, the proof of the next lemma is the same as that of
\cite[Theorem 5.83(b)]{SV}. So we omit the proof.

\begin{lemma}\label{l:VGD}
If
$(y_j)_{j \ge 1}$ is a sequence of points in $D$ converging to $z\in \partial D$,
then for every $x\in D$,
$$
\lim_{j\to \infty} \widetilde V^D \left(\frac{G^D(\cdot, y_j)}{G^D(x_0, y_j)}\right)(x)=
\lim_{j\to \infty} \frac{\widetilde V^D(G^D(\cdot, y_j))(x)}{G^D(x_0, y_j)}
=\widetilde V^D(M^D(\cdot, z))(x)\, .
$$
\end{lemma}

Let us define the function $H^D_Y(x,z):=
\widetilde V^D (M^D(\cdot, z))(x)$ on $D\times \partial D$.
Let $(y_j)$ be a sequence of points in $D$ converging
to $z\in \partial D$, then
from Lemma \ref{l:VGD} we get that
\begin{equation}\label{e5.1}
H^D_Y(x,z)=\lim_{j\to \infty} \frac{\widetilde V^D(G^D(\cdot, y_j))(x)}{G^D(x_0, y_j)}=\lim_{j\to \infty} \frac{U^D(x,y_j)}{G^D(x_0, y_j)}\, ,
\end{equation}
where the last equality follows from Lemma \ref{l:VGD0}.
In particular, there exists the limit
\begin{align}\label{e5.2}
\lim_{j\to \infty} \frac{U^D(x_0,y_j)}{G^D(x_0, y_j)}
= H^D_Y(x_0, z)\, .
\end{align}
Now we define a function $\overline M^D_Y$ on $D\times \partial D$ by
\begin{align}\label{e5.3}
\overline M^D_Y(x, z):=\frac{H^D_Y(x, z)}{H^D_Y(x_0, z)},
\quad x\in D, z\in \partial D.
\end{align}
From the definition above and (\ref{e5.1})--\eqref{e5.2},  we can easily see that
\begin{align}\label{e5.4}
\lim_{D\ni y\to z}\frac{U^D(x,y)}{U^D(x_0, y)}
=\overline M^D_Y(x, z), \quad x\in D, z\in \partial D.
\end{align}
Thus we have proved the following result.
\begin{prop}\label{p:extension-M}
The function $M^D_Y(\cdot, \cdot)$ can be extended  from $D \times D$ to
$D \times \overline{D}$
so that for each $z\in \partial D$ we have that
$$
\overline M^D_Y(x,z)=\lim_{y\to z}M^D_Y(x,y)=\lim_{y\to z}\frac{U^D(x,y)}{U^D(x_0,y)}\, .
$$
\end{prop}

The following two types of sharp two-sided estimates for $\overline M^D_Y(x,z)$ now
follow easily from Theorems \ref{t:gfe} and \ref{t:gfe:nn}.

\begin{thm}\label{t:mke}
Assume that $\phi$ is a  Bernstein function satisfying
{\bf (A1)}--{\bf (A5)}. Let $D\subset \R^d$
be a bounded $C^{1,1}$ domain, or
a $C^{1,1}$ domain with compact complement or
domain above the graph of a bounded $C^{1,1}$
function.
For every $M>0$ and  $z\in \partial D$, there
exists a constant $c=c(M, z)\ge 1$ such that for all $x\in D$ with $|x-z|\le M$,
\begin{align}\label{e:mke}
c^{-1}\frac{\delta(x)\phi'(|x-z|^{-2})}{|x-z|^{d+4}\phi(|x-z|^{-2})^2}\le
\overline  M^D_Y(x,z)\le c \frac{\delta(x)\phi'(|x-z|^{-2})}{|x-z|^{d+4}\phi(|x-z|^{-2})^2}\, .
\end{align}
\end{thm}

Note that the constant $c$ in Theorem \ref{t:mke} will in general depend on
$z\in \partial D$. This is inconsequential, because the
point $z$ will always be fixed.

\begin{thm}\label{t:mke2}
Assume that $\phi$ is a  Bernstein function satisfying
{\bf (A1)}--{\bf (A3)} and {\bf (A6)}.
(1) Let $D\subset \R^d$ be a
domain above the graph of a bounded $C^{1,1}$
function.
There exists a constant $c_1\ge 1$ such that for all $x \in D$ and $z \in \partial D$,
\begin{align}\label{e:mke2}
c_1^{-1}{\delta(x)}
\frac{u(|x-z|^2)|x_0-z|^{d}}
{u(|x_0-z|^2)|x-z|^{d}}\le
\overline M^D_Y(x,z)\le c_1{\delta(x)}
\frac{u(|x-z|^2)|x_0-z|^{d}}
{u(|x_0-z|^2)|x-z|^{d}}\, .
\end{align}
(2) Let $D\subset \R^d$ be a
$C^{1,1}$ domain with compact complement.
There exists a constant $c_2\ge 1$ such that
for all $x \in D$ and $z \in \partial D$,
\begin{align}\label{M:mke4}
&c_2^{-1}\left(\frac{\delta(x)}{|x-z|\wedge 1} \wedge 1 \right)
\left( \frac{|x_0-z|\wedge 1}{|x-z|\wedge 1} \right) \frac{u(|x-z|^2)|x_0-z|^{d-2}}
{u(|x_0-z|^2)|x-z|^{d-2}} \le
\overline M^D_Y(x,z) \nn\\
&\qquad \le c_2 \left(\frac{\delta(x)}{|x-z|\wedge 1} \wedge 1 \right)
\left(  \frac{|x_0-z|\wedge 1}{|x-z|\wedge 1} \right) \frac{u(|x-z|^2)|x_0-z|^{d-2}}
{u(|x_0-z|^2)|x-z|^{d-2}} \, .
\end{align}
 \end{thm}

 \begin{remark}\label{r:mke2} {\rm
 (1) Theorem \ref{t:mke} in particular implies that $\overline M^D_Y(\cdot , z_1)$ differs
from $\overline M^D_Y(\cdot, z_2)$ if $z_1$ and $z_2$ are two different points on $\partial D$.

(2) From  Theorem \ref{t:mke2}, we have $\lim_{D\ni x\to \infty} \overline M^D_Y(x,z)=0$ for any $z\in \partial D$.
In fact, for $|x-z|\ge |z-x_0|$ we have
${u(|x-z|)}\le {u(|x_0-z|)} $.
It is clear that
$$
\limsup_{D\ni x\to \infty} \left({\delta(x)}
\frac{|x_0-z|^{2}}
{|x-z|^{2}}+\frac{|x_0-z|\wedge 1}{|x-z|\wedge 1}\right) \le \limsup_{D\ni x\to \infty} \left(
\frac{|x_0-z|^2}
{|x-z|}+\frac{|x_0-z|\wedge 1}{|x-z|\wedge 1}\right)<\infty \, .
$$
Thus, in both cases,
\begin{align}\label{M:mke5}
\limsup_{D\ni x\to \infty} \overline M^D_Y(x,z)\le c \limsup_{D\ni x\to \infty}
\frac{u(|x-z|^2)|x_0-z|^{d-2}}
{u(|x_0-z|^2)|x-z|^{d-2}}\le c \limsup_{D\ni x\to \infty} \frac{|x_0-z|^{d-2}}
{|x-z|^{d-2}}=0.
\end{align}
}
\end{remark}

Using  the continuity of $U^D$ in the extended sense (Proposition \ref{p:U-continuous})
and the upper bound
in \eqref{e:free-gfe-reverse}, one can check that $Y^D$ satisfies Hypothesis (B) in \cite{KW}.
Therefore, $D$  has a Martin boundary $\partial_M D$ with respect to $Y^D$ satisfying the following properties:
\begin{description}
\item{\bf (M1)} $D\cup \partial_M D$ is
a compact metric space (with the metric denoted by $d$);
\item{\bf (M2)} $D$ is open and dense in $D\cup \partial_M D$,  and its relative topology coincides with its original topology;
\item{\bf (M3)}  $M^D_Y(x ,\, \cdot\,)$ can be uniquely extended  to $\partial_M D$ in such a way that
\begin{description}
\item{(a)}
$ M^D_Y(x, y) $ converges to $M^D_Y(x, w)$ as $y\to w \in \partial_M D$ in the Martin topology;
\item{(b)} for each $ w \in D\cup \partial_M D$ the function $x \to M^D_Y(x, w)$  is excessive with respect to $Y^D$;
\item{(c)} the function $(x,w) \to M^D_Y(x, w)$ is jointly continuous on $D\times ((D\setminus\{x_0\})\cup \partial_M D)$ in the Martin topology and
\item{(d)} $M^D_Y(\cdot,w_1)\not=M^D_Y(\cdot, w_2)$ if $w_1 \not= w_2$ and $w_1, w_2 \in \partial_M D$.
\end{description}
\end{description}

Recall that  a positive harmonic function $f$ for  $Y^{D}$ is
minimal if, whenever $h$ is a positive harmonic function for
$Y^{D}$ with $h\le f$ on $D$, one must have $f=ch$ for some
constant $c$. A point $z\in \partial_M D$ is called a minimal
Martin boundary point
if $M^D_Y(\cdot, z)$ is a minimal harmonic function for $Y^D$.
The minimal Martin boundary of $Y^D$ is denoted by $\partial_m D$.

We will say that a point $w\in \partial_M D$ is a finite Martin
boundary point if there exists a bounded sequence
$(y_n)_{n\ge 1}\subset D$ converging to $w$ in the Martin topology.
Recall that a point $w$ on the Martin boundary $\partial_MD$ of $D$ is said to be associated
with $z\in \partial D$ if there is a sequence $(y_n)_{n\ge 1}\subset D$ converging to $w$
in the Martin topology and to $z$ in the Euclidean topology. The set of Martin
boundary points associated with $z$ is denoted by $\partial_M^zD$.

By using Proposition \ref{p:extension-M}, the proof of next lemma is same  as that  of \cite[Lemma 3.6]{KSV10}. Thus we omit it.

\begin{prop}\label{p:infinite-mb}
For any $z\in \partial D$, $\partial_M^zD$ consists of exactly one
point $w$ and $M^D_Y(\cdot, w)=\overline M^D_Y(\cdot, z)$.
\end{prop}

Because of the proposition above, we will also use $z$ to denote the point on the
Martin boundary $\partial_M^zD$ associated with $z\in \partial D $.
Note that it follows from the proof of \cite[Lemmas 3.6]{KSV10} that if $(y_n)_{n\ge 1}$
converges to $z\in \partial D$ in the Euclidean topology, then it also converges to $z$ in the Martin topology.

In the remainder of this section,  we fix  $z\in \partial D$.
The proof of the next result is same  as that of \cite[Lemma 3.8]{KSV10}.
Thus we omit the proof.

\begin{lemma}\label{l:mk-integrability}
For every bounded open $O\subset \overline{O}\subset D$ and every
$x\in D$, $M^D_Y(Y^D_{\tau_O}, z)$ is $\P_x$-integrable.
\end{lemma}

Using the results above, we can get the following result.

\begin{lemma}\label{l:mk-harmonic}
Suppose that $\phi$ is a  Bernstein function satisfying
{\bf (A1)}--{\bf (A6)}.
For any $x\in D$ and $r \in (0, R\wedge(\delta(x)/2)]$,
$$
    M^D_Y(x, z)=\E_x[M^D_Y(Y^D_{\tau_{B(x, r)}}, z)]\, .
$$
\end{lemma}

\pf
Recall that  $D$ satisfies the interior and exterior balls conditions with radius
$R$. Thus,
for all $r\in(0,R]$,
there is a ball  $B(A_r(z) , r/2)  \subset D\cap B(z,r)$.
 Fix $x\in D$ and a positive $r<R\wedge \frac{\delta(x)}2$. Let
$$
\eta_m:=2^{-2m}r \quad \mbox{and  }\ z_m=A_{\eta_m}(z),
\quad m=0, 1, \dots.
$$
Note that
$$
B(z_m, \eta_{m+1})\subset D\cap B(z, 2^{-1}\eta_m)\subset D\cap B(z,\eta_m)\subset
D\cap B(z, r)\subset D\setminus B(x, r)
$$
for all $m\ge 0$. Thus by the harmonicity of $M^D_Y(\cdot, z_m)$, we have
$$
M^D_Y(x, z_m)=\E_x\left[M^D_Y(Y^D_{\tau_{B(x, r)}}, z_m)\right].
$$
Choose $m_0=m_0 \ge 2$
such that $\eta_{m_0} <\delta(x_0)/4$.

To prove the lemma, it suffices to show that
$\{M^D_Y(Y^D_{\tau_{B(x, r)}}, z_m): m\ge m_0\}$ is $\P_x$-uniformly
integrable. Fix an arbitrary $\eps>0$.
We first note that  if $D$ is  unbounded,
by   Theorem \ref{t:gfe:nn}
there exists $L \ge 2r\vee 2$ such that for every $m \ge m_0$ and $w\in D\setminus B(z, L)$,
\begin{align*}
&\frac{U^D(w, z_m)}
 {U^D(x_0, z_m)} \le
 \frac{c}{\delta(z_m)}
 \left(\frac{\delta(z_m)}{|w-z_m|\wedge 1} \wedge 1 \right)
\left(\frac{\delta(w)}{|w-z_m|\wedge 1} \wedge 1 \right)\frac{u(|w-z_m|^2)}
{|w-z_m|^{d-2}}\\
&\le
 \frac{c}{\delta(z_m)}
 \left(\delta(z_m) \wedge 1 \right)
\left(\delta(w) \wedge 1 \right)\frac{u(|w-z_m|^2)}
{|w-z_m|^{d-2}}\le
c\frac{u(|w-z_m|^2)}
{|w-z_m|^{d-2}}\\
&\le
c \frac{\phi'(|w-z_m|^{-2})}
{|w-z_m|^{d+2}\phi(|w-z_m|^{-2})^2} \le
c \frac{\phi'((L/2)^{-2})}
{(L/2)^{d+2}\phi((L/2)^{-2})^2}
\le \frac{\eps}{4}.
\end{align*}
In the above inequalities, we have used Lemma \ref{l:properties-of-bf}(b).
If $D$ is a bounded domain we simply take
$L=2\mathrm{diam}(D)$ so that $D\setminus B(z, L)= \emptyset$.
   Thus
\begin{align}\label{e:KSV15.8-1}
&\E_x\left[M^D_Y(Y^D_{\tau_{B(x, r)}}, z_m);
Y^D_{\tau_{B(x, r)}}\in D\setminus B(z, L)\right] \le \frac{\eps}{4}\, .
\end{align}
By Theorem \ref{t:gfe},
there exist $m_1 \ge m_0$ and $c_1=c_1(L)>0$
such that   for every $w\in (D \cap B(z, L)) \setminus  B(z, \eta_m)$ and $y\in D\cap B(z, \eta_{m+1})$,
$$
M^D_Y(w, z_m)\le c_1M^D_Y(w, y), \qquad m\ge m_1.
$$
Letting $y\to z$ we get
\begin{align}\label{e:KSV15.8}
M^D_Y(w, z_m)\le c_1M^D_Y(w, z), \qquad m\ge m_1, w\in (D \cap B(z, L)) \setminus  B(z, \eta_m).
\end{align}
Since $M^D_Y(Y^D_{\tau_{B(x, r)}}, z)$ is $\P_x$-integrable by
Lemma \ref{l:mk-integrability},  there is an $N_0=N_0(\eps)>1$ such that
\begin{align}\label{e:KSV15.9}
\E_x\left[M^D_Y(Y^D_{\tau_{B(x, r)}}, z); M^D_Y(Y^D_{\tau_{B(x, r)}}, z)>N_0/c_1\right]<\frac{\epsilon}{2c_1}.
\end{align}
 By \eqref{e:KSV15.8-1},  \eqref{e:KSV15.8} and \eqref{e:KSV15.9},
\begin{eqnarray*}
&&\E_x\left[M^D_Y(Y^D_{\tau_{B(x, r)}}, z_m); M^D_Y(Y^D_{\tau_{B(x, r)}}, z_m)>N_0 \mbox{ and }
 Y^D_{\tau_{B(x, r)}}\in D\setminus B(z, \eta_m)\right]\\
 &\le&\E_x\left[M^D_Y(Y^D_{\tau_{B(x, r)}}, z_m); M^D_Y(Y^D_{\tau_{B(x, r)}}, z_m)>N_0 \mbox{ and }
 Y^D_{\tau_{B(x, r)}}\in (D \cap B(z, L)) \setminus  B(z, \eta_m)\right]\\
 &&+\E_x\left[M^D_Y(Y^D_{\tau_{B(x, r)}}, z_m);
 Y^D_{\tau_{B(x, r)}}\in D\setminus B(z, L)\right]\\
&\le &c_1\E_x\left[M^D_Y(Y^D_{\tau_{B(x, r)}}, z); c_1M^D_Y(Y^D_{\tau_{B(x, r)}}, z)>N_0\right] + \frac{\eps}{4}
<c_1\frac{\epsilon}{2c_1}+ \frac{\eps}{4}=\frac{3\epsilon}{4}.
\end{eqnarray*}

By \eqref{e:KSV7-4.7}, we have for $m\ge m_1$,
\begin{eqnarray*}
&&\E_x\left[M^D_Y(Y^D_{\tau_{B(x, r)}}, z_m);Y^D_{\tau_{B(x, r)}}\in D\cap B(z, \eta_m)\right]\\
&& =\int_{D\cap B(z, \eta_m)}M^D_Y(w, z_m)K^{D, B(x, r)}(x, w)dw\\
&&\le c_2 \phi(r^{-2})^{-1}\int_{D\cap B(z, \eta_m)}M^D_Y(w, z_m)
 \delta(w)\frac{\phi'((|w-x|-r)^{-2})}
{(|w-x|-r)^{d+3}} dw.
\end{eqnarray*}
Since $|w-x|\ge |x-z|-|z-w|\ge \delta(x)-\eta_m\ge
\frac{7}4r$, applying Lemma \ref{l:properties-of-bf}(a)--(c), we get that
\begin{eqnarray}
&&\E_x\left[M^D_Y(Y^D_{\tau_{B(x, r)}}, z_m);Y^D_{\tau_{B(x, r)}}\in D\cap
B(z, \eta_m)\right]\nonumber\\
&&\le c_3 r^{-d-3}\phi'(((3r/4)^{-2})\phi((3r/4)^{-2})^{-1}
\int_{D\cap B(z, \eta_m)}M^D_Y(w, z_m)\delta(w)dw\nonumber\\
&&\le c_4r^{-d-3}\phi'(r^{-2})\phi(r^{-2})^{-1}U^D(x_0, z_m)^{-1}\int_{D\cap B(z, \eta_m)}U^D(w, z_m)\delta(w)dw.
\label{e:KSV15.10}
\end{eqnarray}
Note that, by Theorem \ref{t:gfe} ,
\begin{align}\label{e:KSV15.11}
U^D(x_0, z_m)^{-1}\le \frac{c_5}{\eta_{m}}
\end{align}
and by \eqref{e:free-gfe-reverse}
\begin{eqnarray}
&&\int_{D\cap B(z, \eta_m)}\delta(w)U^D(w, z_m)dw \le
\int_{D\cap B(z, \eta_m)}\delta(w)G_X(w, z_m)dw
\nn\\
&\le& c_6 \eta_m \int_{D\cap B(z, \eta_m)}
\frac{\phi'(|w-z_m|^{-2})}
{|w-z_m|^{d+2}\phi(|w-z_m|^{-2})^2}dw\nonumber\\
&\le&
c_6  \eta_m\int_{B(z_m, 2\eta_m)}
\frac{\phi'(|w-z_m|^{-2})}
{|w-z_m|^{d+2}\phi(|w-z_m|^{-2})^2}dw\nonumber\\
&=& c_6\eta_m  \int_{B(0, 2\eta_m)}\frac{\phi'(|w|^{-2})}
{|w|^{d+2}\phi(|w|^{-2})^2}dw
= c_7\eta_m\int_0^{2\eta_m}\frac{\phi'(r^{-2})}
{r^3\phi(r^{-2})^2}dr
\nonumber\\
&=&c_7\eta_m\int_0^{2\eta_m}\frac{d}{dr}\left(\frac{1}{\phi(r^{-2})}\right)dr
\le c_8  \eta_m \phi((2\eta_m)^{-2})^{-1}.\label{e:KSV15.12}
\end{eqnarray}
It follows from \eqref{e:KSV15.10}--\eqref{e:KSV15.12} that
\begin{eqnarray*}
&&\E_x\left[M^D_Y(Y^D_{\tau_{B(x, r)}}, z_m);Y^D_{\tau_{B(x, r)}}\in D\cap B(z, \eta_m)\right]\\
&&\le c_9r^{-d-3}\phi'(r^{-2}) \phi(r^{-2})^{-1}\frac{1}
{\phi((2\eta_m)^{-2})} \le \frac{c(r)}
{\phi((2\eta_m)^{-2})} .
\end{eqnarray*}
Thus there exists $m_2 \ge m_1$ such that for all $m\ge m_2$,
$$
\E_x\left[M^D_Y(Y^D_{\tau_{B(x, r)}}, z_m);Y^D_{\tau_{B(x, r)}}\in D\cap B(z, \eta_m)\right]\le \frac{\epsilon}4.
$$
Consequently, for all $m\ge m_2$,
$$
\E_x\left[M^D_Y(Y^D_{\tau_{B(x, r)}}, z_m);M^D(Y^D_{\tau_{B(x, r)}}, z_m)>N\right]\le \epsilon,
$$
which implies that $\{M^D_Y(Y^D_{\tau_{B(x, r)}}, z_m): m\ge m_0\}$ is $\P_x$-uniformly
integrable.
\qed

Using this, we can easily get the following

\begin{thm}\label{t:mk-harmonic}
Suppose that $\phi$ is a  Bernstein function satisfying
{\bf (A1)}--{\bf (A6)}.
The function $M^D(\cdot, z)$ is harmonic in $D$ with respect to $Y^D$.
\end{thm}

\pf The proof is the same as that of \cite[Theorem 3.10]{KSV10}. \qed

\begin{thm}\label{t:m}
Suppose that $\phi$ is a  Bernstein function satisfying
{\bf (A1)}--{\bf (A6)}.
Every point $z$ on $\partial D$ is a minimal Martin boundary
point.
\end{thm}

\pf
Fix $z \in \partial D$ and let $h$ be a positive harmonic function for  $Y^D$
such that $h\le M^D_Y(\cdot, z)$. By
the Martin representation in \cite{KW},
there is a finite measure on $\partial_M D$ such that
$$
    h(x)=\int_{\partial_M D}M^D_Y(x,w)\, \mu(dw)=\int_{\partial_M D\setminus\{z\}}M^D_Y(x,w)\, \mu(dw)+
    M^D_Y(x, z)\mu(\{z\})\, .
$$
In particular, $\mu(\partial_M D)=h(x_0)\le
M^D_Y(x_0, z)=1$ (because of the normalization at $x_0$). Hence, $\mu$ is a sub-probability measure.

For $\epsilon >0$, put
$K_{\epsilon}:=\left\{w\in \partial_M D: d(w, z) \ge \epsilon\right\}$.
Then $K_{\epsilon}$ is a compact subset of $\partial_M D$.
Define
\begin{align}\label{d:definition-u}
    u(x):=\int_{ K_{\epsilon} }M^D_Y(x,w)\, \mu(dw).
\end{align}
Then $u$ is a positive harmonic function with respect to  $Y^{D}$
satisfying
\begin{align}\label{e:newm1}
u(x)\le h(x)-\mu(\{z\})M^D_Y(x,
z)\le \big(1-\mu(\{z\})\big)M^D_Y(x, z)\, .
\end{align}
By (M3)(c), our estimates in
Theorems \ref{t:mke} and \ref{t:mke2}
and the fact $\lim_{D\ni x\to \infty} M^D_Y(x,z)=0$ (cf.~Remark \ref{r:mke2})
we see from \eqref{d:definition-u} and \eqref{e:newm1} that
$u$ is bounded,   $\lim_{D\ni x\to w}u(x)=0$ for every  $w\in \partial D$ and
$\lim_{D\ni x\to \infty}u(x)=0$.
Therefore by the harmonicity of $u$,  $u\equiv 0$ in $D$.

We see from \eqref{d:definition-u} that $\nu=\mu_{| K_{\epsilon}}=0$. Since
$\epsilon >0$ was arbitrary and $\partial_M D\setminus\{z\}=\cup_{\epsilon >0} K_{\epsilon}$,
we see that $\mu_{|\partial_M D\setminus\{z\}}=0$. Hence $h=\mu(\{z\})
M^D_Y(\cdot, z)$ showing that $
M^D_Y(\cdot, z)$ is minimal.
\qed

Combining Remark \ref{r:mke2}(1) and Theorem \ref{t:m}, we conclude that
\begin{thm}\label{mbid}
Suppose that $\phi$ is a  Bernstein function satisfying
{\bf (A1)}--{\bf (A6)}.
The finite part of the minimal Martin boundary of $D$ and the finite
part of the  Martin boundary of $D$ both coincide with  the Euclidean boundary $\partial D$ of $D$.
\end{thm}

We conclude this section with following inequality, which will be used in Section \ref{s:mtf}.

\begin{corollary}\label{c:U-M}
Fix $z \in \partial D$ and assume that
$x_0\in D\cap B(z, R)$ satisfies $R/4<\delta(x_0)<R$ and $M^D_Y$
is the Martin kernel of $D$ based on $x_0$. Then
there exists $c=c(z)>0$ such that for all $x,y\in B(z, R/4)$ with
$\frac34 |x-z|\le |x-y|$,
\begin{align}\label{e:U-M}
\frac{U^D(x,y)}{M^D_Y(x,z)}\,\le \,c\, U^D(x_0,y)\, .
\end{align}
\end{corollary}
\pf It follows from Theorem \ref{t:gfe} and Theorem \ref{t:mke} that
\begin{eqnarray*}
U^D(x,y)&\asymp & \delta(x)\delta(y)|x-y|^{-d-4}\phi'(|x-y|^{-2})\phi(|x-y|^{-2})^{-2}\, ,\\
M^D_Y(x,z)& \asymp & \delta(x) |x-z|^{-d-4}\phi'(|x-z|^{-2})\phi(|x-z|^{-2})^{-2}\, ,\\
U^D(x_0,y)&\asymp & \delta(y)|x_0-y|^{-d-4}\phi'(|x_0-y|^{-2})\phi(|x_0-y|^{-2})^{-2}
\,\asymp \, \delta(y)\, .
\end{eqnarray*}
Since $|x_0-y|\ge R/4$ and $r\mapsto r^{-d-4}\phi'(r^{-2})\phi(r^{-2})^{-2}$ is decreasing,
we can estimate $U^D(x_0,y)
\ge c_1 \delta(y)$.
Using the monotonicity of $r\mapsto r^{-d-4}\phi'(r^{-2})\phi(r^{-2})^{-2}$, we get
$$
\frac{\phi'(|x-y|^{-2})}{|x-y|^{d+4}\phi(|x-y|^{-2})^{2}}\le c
\, \frac{\phi'(((3|x-z|)/4)^{-2})}{((3|x-z|)/4)^{d+4}
\phi(((3|x-z|)/4)^{-2})^{2}}\, .
$$
Applying Lemma \ref{l:properties-of-bf}(c)
we get that
${U^D(x,y)}/{M^D_Y(x,z)}\le c_1 \delta(y)$.
This
completes the proof. \qed

\section{Quasi-additivity of capacity}\label{s:quasi}

Throughout this section we assume that $\phi$ is a
Bernstein function satisfying {\bf (A1)}--{\bf (A5)}.
Let $\mathrm{Cap}$ denote the capacity with respect to the
subordinate Brownian motion $X$ and $\mathrm{Cap}_D$ the capacity
with respect to the subordinate killed Brownian motion $Y^D$.
The goal of this section is to prove that $\mathrm{Cap}_D$ is
quasi-additive with respect to Whitney decompositions of $D$.

We start with the following  inequality: There exist
positive constants $c_1<c_2$ such that
\begin{align}\label{e:capacity-estimate}
c_1 r^d \phi(r^{-2})\le \mathrm{Cap}(\overline{B(0,r)})
\le c_2 r^d \phi(r^{-2})\, ,\quad \text{for every }r\in (0,1]\, .
\end{align}
Using
\eqref{e:free-gfe-reverse}, the proof of \eqref{e:capacity-estimate} is the same as that of  \cite[Proposition 5.2]{KSV11}. Thus we omit the proof.

For any open set $D\subset \R^d$, let $\SS(D)$  denote the
collection of all excessive functions with respect to $Y^D$ and
let $\SS^c(D)$ be the family of positive  functions in $\SS(D)$
which are continuous in the extended sense.
For any $v\in \SS(D)$ and $E\subset D$, the reduced function of $v$
relative to $E$ in $D$ is defined by
\begin{align}\label{e:RE}
R^E_v(x)=\inf\{w(x): w\in \SS(D) \mbox{ and } w\ge v
\mbox{ on } E\}, \qquad x\in \R^d.
\end{align}
The lower semi-continuous regularization $\wh R^E_v$ of
$R^E_v$ is called the balayage of $v$ relative to $E$ in $D$.
Note that the killed Brownian motion $W^D$ is a strongly Feller process.
Thus it follows by \cite[Proposition V.3.3]{BH} that
the semigroup of $Y^D$ also has strong Feller property.
So it follows easily
from \cite[Proposition V.2.2]{BH} that the cone of excessive
functions $\SS(D)$ is a balayage space in the sense of  \cite{BH}.

In the remainder of  this section we assume that
$D\subset \R^d$ is either a bounded $C^{1,1}$ domain, or
a $C^{1,1}$ domain with compact complement or a
domain above the graph of a bounded  $C^{1,1}$ function.

Given $v\in \SS^c(D)$, define a kernel $k_v:D\times D\to [0,\infty]$ by
\begin{align}\label{e:ku}
k_v(x,y):=\frac{U^D(x,y)}{v(x)v(y)}\, ,\qquad x,y\in D\,.
\end{align}
We will later consider $v(y)=U^D(y,x_0)\wedge 1$.
Note that $k_v(x,y)$ is jointly lower semi-continuous on $D \times D$
by the joint lower semi-continuity of $U^D$, cf.~Proposition \ref{p:U-continuous}, and the assumptions
that $v$ is positive and continuous in the extended sense.
For a measure $\lambda$ on $D$ let $\lambda_v(dy):=\lambda(dy)/v(y)$. Then
$$
k_v\lambda (x):=\int_D k_v(x,y)\, \lambda(dy)=\int_D \frac{U^D(x,y)}{v(x)v(y)}\,
\lambda(dy)=\frac{1}{v(x)}\int_D U^D(x,y)\, \frac{\lambda(dy)}{v(y)}
=\frac{1}{v(x)}U^D \lambda_v(dy)\, .
$$
We define a capacity with respect to the kernel $k_v$ as follows:
\begin{eqnarray*}
\sC_v(E):=\inf\{\|\lambda\|: k_v\lambda \ge 1 \textrm{ on }E\}\, ,
\qquad E\subset D\, ,
\end{eqnarray*}
where $\|\lambda \|$ denotes the total mass of the measure
$\lambda$ on $D$. The following dual representation of the
capacity of compact sets can be found in \cite[Th\'eor\`eme 1.1]{Fug}:
\begin{align}\label{e:fuglede-equality}
\sC_v(K)=\sup\{\mu(K):\, \mu(D\setminus K)=0, k_v\mu \le 1\ \textrm{on }D\}\, .
\end{align}
For a compact set $K\subset D$, consider  the balayage $\wh{R}_v^K$.
Being a potential, $\wh{R}_v^K=U^D \lambda^{K, v}$ for a measure
$\lambda^{K, v}$ supported in $K$.
Recall that $(\EE^D, \DD(\EE^D))$ is  the Dirichlet form associated with $Y^D$.
Define the Green energy of $K$ (with respect to $v$) by
$$
\gamma_v(K):=\int_D \int_D U^D(x,y) \lambda^{K, v}(dx)\, \lambda^{K, v}(dy) =
\int_D U^D \lambda^{K, v}(x)\, \lambda^{K, v}(dx)=\EE^D(U^D \lambda^{K, v},U^D \lambda^{K, v})\, .
$$
As usual, this definition of energy is extended first to open and then
to Borel subsets of $D$. By following the proof of \cite[Proposition 5.3]{KSV11}
we see that for all Borel subsets $E\subset D$ it holds that
\begin{align}\label{e:gamma=c}
\gamma_v(E)=\sC_v(E)\, .
\end{align}
Note that in case $v\equiv 1$, $\gamma_1(E)=\sC_1(E)=\cp_D(E)$.

Let $\{Q_j\}_{j\ge 1}$ be a Whitney decomposition of $D$.
Recall that $x_j$ is the center of $Q_j$ and $Q_j^*$ the
interior of the double of $Q_j$. Then $\{Q_j, Q_j^*\}$ is
a quasi-disjoint decomposition of $D$ in the sense of \cite[pp.~146-147]{AE}.
\begin{defn}\label{kernelHI}
A kernel $k:D\times D\to [0,+\infty]$ is said to satisfy the
local Harnack property with localization constant $r_1>0$
with respect to $\{Q_j, Q_j^*\}$ if
\begin{align}\label{e:hpk}
k(x,y)\asymp k(x',y)\, , \textrm{ for all }x,x'\in Q_j \textrm{  and all }y\in D\setminus Q_j^*\, ,
\end{align}
for all cubes $Q_j$ of diameter less than $r_1$.
\end{defn}

\begin{defn}\label{d:sihi}
A function $v:D\to (0,\infty)$ is said to satisfy the local scale invariant
Harnack inequality with localization constant $r_1>0$
with respect to $\{Q_j\}$
if there exists $c>0$ such that
\begin{align}\label{e:u-scale-inv-har}
\sup_{Q_j}v\le c \inf_{Q_j}v\, ,\quad \textrm{for all } Q_j \text{ with }\mathrm{diam}(Q_j)<r_1\, .
\end{align}
\end{defn}

\begin{lemma}\label{l:sihi-khp}
If $v\in \SS^c(D)$ satisfies the local scale invariant Harnack inequality
with localization constant $r_1>0$ with respect to $\{Q_j\}$,
then the kernel $k_v$ satisfies the local Harnack property
with localization constant $r_1>0$ with respect to $\{Q_j, Q_j^*\}$.

\end{lemma}
\pf
This is an immediate consequence of Corollary \ref{c:whitney-green}(i).
\qed

Typical examples of positive continuous excessive functions $v$
that satisfy the scale invariant Harnack inequality
are functions $v\equiv 1$ and $v=U^D(\cdot,x_0)\wedge c$ with $x_0\in D$  and $c>0$ fixed.

\begin{lemma}\label{l:cap-comparability}
For every $M>0$, there exists  a constant $c=c(M)\in (0,1)$ such that
\begin{align}\label{e:cap-comparability}
c\, \cp_D(Q_j)\le \cp(Q_j)\le \cp_D(Q_j)
\end{align}for all Whitney cubes whose diameter is less than $M$.
\end{lemma}

\pf By
\eqref{e:fuglede-equality} and \eqref{e:gamma=c} we have that
for every compact set $K\subset D$,
\begin{eqnarray*}
\cp_D(K)=\sup\{\mu(K):\, \mathrm{supp}(\mu)\subset K, U^D\mu
\le 1 \text{ on }D\}\, .
\end{eqnarray*}
If $\mathrm{supp}(\mu)\subset K$ and $G_X\mu\le 1$ on $\R^d$,
then clearly $U^D \mu\le 1$ on $D$.
This implies that $\cp(K)\le \cp_D(K)$ for all compact subset $K\subset D$,
in particular for each Whitney cube $Q_j$.

Let $\mu$ be the capacitary measure of $Q_j$ (with respect to $Y^D$), i.e., $\mu(Q_j)=\cp_D(Q_j)$ and $U^D\mu \le 1$.
Then  by Corollary \ref{c:whitney-green}(ii) for every $x\in Q_j$ we have
$$
1\ge  U^D \mu(x)=\int_{Q_j}U^D(x,y)\, \mu(dy)\ge  \int_{Q_j}c G_X(x,y)\, \mu(dy)=G_X(c\mu)(x)\, .
$$
By the maximum principle it follows that $G_X(c\mu)\le 1$ everywhere on $\R^d$. Hence,
$\cp(Q_j)\ge (c\mu)(Q_j)=c\, \cp_D(Q_j)$.
\qed

\begin{lemma}\label{l:comp-gamma-cap}
Suppose that $v\in \SS^c(D)$ is a function satisfying the local scale invariant
Harnack inequality with localization constant $r_1>0$ with respect to $Y^D$.
Then for every $Q_j$ of diameter less than $r_1$ and every
$E\subset Q_j$ it holds that
\begin{align}\label{e:comp-gamma-cap}
\gamma_v(E)\asymp v(x_j)^2 \cp_D(E)\, .
\end{align}
\end{lemma}
\pf The proof is same as the proof of \cite[Lemma 5.8(i)]{KSV11}. \qed

\begin{defn}\label{d:comp_measure}
Let  $\{Q_j\}$ be  a Whitney decomposition of $D$ and $v\in \SS^c(D)$. A Borel measure $\sigma$ on $D$
is locally comparable to the capacity $\sC_v$ with respect to $\{Q_j\}$ at $z \in \partial D$
if there exists $r, c>0$ such that
\begin{eqnarray*}
& &\sigma(Q_j)\asymp \sC_v(Q_j),\quad \textrm{for all }Q_j  \textrm{ with } Q_j \cap B(z, r)
\not= \emptyset\, ,\\
& &\sigma(E)\le c\,  \sC_v(E),\quad \textrm{for all Borel }E \subset D \cap B(z, 2 r).
\end{eqnarray*}
\end{defn}

Recall that
$(\EE^D, \DD(\EE^D))$ is  the Dirichlet form associated with $Y^D$.

\begin{lemma}\label{l:hardy} (Local Hardy's inequality)
There exist constants $c>0$ and
$r >0$ such that
for every $v\in \DD(\EE^D)$ and $z \in \partial D$,
\begin{align}\label{e:hardy}
\EE^D(v,v)\ge c \int_{D \cap B(z, r)} v(x)^2 \phi(\delta(x)^{-2})\, dx\, .
\end{align}
\end{lemma}
\pf
Since $D$ is a $C^{1,1}$ domain, there exist $b_1>1$,  $R_1>0$ and a cone $C$ whose vertex is at the origin, such that
for every $z \in \partial D$ and $x\in D \cap B(z, b_1R_1/2)$, there exists $\wh C$, which is a rotation of $C$, such that
 \begin{align}
 \label{e:he11}
  (\wh C+x) \cap \{b_1\delta (x)< |x-y|<R_1\}\subset   D^c.
 \end{align}
 Choose $r \in (0, b_1R_1/2)$ small that $\phi((b_1r)^{-2})  \ge 2 \phi(R_1^{-2})$.

Fix $v\in \DD(\EE^D)$ and $z \in \partial D$.
By \eqref{e:dirichlet-form-Y} and \eqref{e:killing-functions},
$$
\EE^D(v,v) \ge \int_{D \cap B(z, r)} v(x)^2 \kappa_D(x)\, dx \ge  \int_{D \cap B(z, r)} v(x)^2 \kappa_D^X(x)\, dx  \, .
$$
Let $x\in D \cap B(z, r)$.
By \eqref{e:kappa-X}, \eqref{e:he11}, and the lower bound in
\eqref{e:levy-density-estimate},
\begin{eqnarray*}
\kappa_D^X(x)&= & \int_{D^c}j(x-y)\,
dy \ge \int_{(\wh C+x) \cap \{b_1\delta (x)< |x-y|<R_1\}}j(x-y)dy\\
& \ge& c_1 \int_{(\wh C+x) \cap \{b_1\delta (x)< |x-y|<R_1\}} |x-y|^{-d-2}\phi'(|x-y|^{-2})
\, dy\\
&\ge & c_2  \int_{b_1 \delta(x)}^{R_1} -\frac{d}{ds}(\phi(s^{-2})) ds
=c_2 (\phi((b_1\delta(x))^{-2})-\phi(R_1^{-2}))\\
&=&2^{-1}c_2 \phi((b_1\delta(x))^{-2}) \ge c_3  \phi(\delta(x)^{-2}) \, .
\end{eqnarray*}
In the second to last inequality we used
$\phi((b_1\delta(x))^{-2}) \ge  \phi((b_1r)^{-2}) \ge 2 \phi(R_1^{-2})$
and, in the last inequality we used \eqref{e:wlsc-substitute}.
 \qed

For $v\in \SS^c(D)$, define
$$
\sigma_v(E):=\int_E v(x)^2 \phi(\delta(x)^{-2})\, dx\, ,\qquad E\subset D\, .
$$

\begin{prop}\label{p:comparable-measure}
Let $v\in \SS^c(D)$ satisfy the local  scale invariant
Harnack inequality with localization constant $r_1>0$ with
respect to the Whitney decomposition $\{Q_j\}$. Then $\sigma_v$ is
 locally comparable to the capacity  $\sC_v$ with respect to $\{Q_j\}$ for every $z \in D$.
\end{prop}
\pf
Fix $z\in\partial D$ and let $\wt r=(r_1 \wedge r_2)/2$ where $r_2$ is the constant $r$ in Lemma \ref{l:hardy}.
Since $v$ satisfies the  local scale invariant
Harnack inequality with localization constant $r_1$,
we have $v\asymp v(x_j)$ on any $Q_j$ of diameter less than $r_1$.
By Lemma \ref{l:comp-gamma-cap},
$\gamma_v(Q_j)\asymp v(x_j)^2 \cp_D(Q_j)$
for any $Q_j$ of diameter less than $\wt r$.
On the other hand,
by Lemma \ref{l:cap-comparability} and \eqref{e:capacity-estimate},
$$
\sigma_v(Q_j)=\int_{Q_j}v(x)^2 \phi(\delta(x)^{-2})\, dx\asymp
v(x_j)^2 \phi\big((\mathrm{diam}(Q_j)\big)^{-2}|Q_j|\asymp \cp(D)
\asymp\cp_D(Q_j)
$$
for all $Q_j$ with $Q_j\cap B(z,\wt r)\neq \emptyset$.
Thus $\gamma_v(Q_j)\asymp \cp_D(Q_j)$.

Using local Hardy's inequality,  Lemma \ref{l:hardy}, for any
Borel subset $E\subset D$ and compact $K\subset E\cap B(z, 2 \wt r)$,
\begin{eqnarray*}
\gamma_v(E)&\ge & \gamma_v(K)=\EE^D(U^D \lambda^{K, v}, U^D \lambda^{K, v})
\ge c_1 \int_K (U^D\lambda^{K, v})(x)^2 \phi(\delta(x)^{-2})\, dx \\
&=& c_1
\int_K v(x)^2 \phi(\delta(x)^{-2})\, dx =
c_1 \sigma_v(K)\, .
\end{eqnarray*}
This proves that $\gamma_v(E)\ge c_1  \sigma_v(E)$. \qed

Now we can repeat the argument in the proof of \cite[Theorem 7.1.3]{AE}
and conclude that $\gamma_v=\sC_v$ is quasi-additive with respect to $\{Q_j\}$.

\begin{prop}\label{p:quasi-additivity}
For any Whitney decomposition $\{Q_j\}$ of $D$ and any $v\in \SS^c(D)$
satisfying the local scale invariant Harnack inequality with respect to
$\{Q_j\}$, the Green energy $\gamma_v$ is  locally quasi-additive with respect to $\{Q_j\}$ for every $z \in \partial D$:
There exist $r, c>0$ such that for every $z \in \partial D$,
$$
c^{-1} \sum_{j\ge 1} \gamma_v(E\cap Q_j)   \le  \gamma_v(E)\le c \sum_{j\ge 1} \gamma_v(E\cap Q_j)\, \quad
\textrm{for all Borel }E\subset D \cap B(z, r).
$$
\end{prop}

\section{Minimal thinness}\label{s:mtf}

Throughout this section, we assume  that $\phi$ is a  Bernstein function satisfying
{\bf (A1)}--{\bf (A6)} and that
$D\subset \R^d$ is either a bounded $C^{1,1}$ domain, or
a $C^{1,1}$ domain with compact complement or
a domain above the graph of a bounded $C^{1,1}$
function. We assume that the $C^{1,1}$ characteristics
of $D$ is $(R, \Lambda)$.

We start this section by recalling the definition of minimal
thinness of a set at a minimal Martin boundary point with
respect to the subordinate killed Brownian motion $Y^D$.

\begin{defn}\label{def:cMthin}
Let $D$ be an open set in  $\R^d$. A set $E\subset D$ is said
to be minimally thin in $D$ at $z\in  \partial_m D$  with
respect to $Y^D$ if $\wh R^E_{M_Y^D(\cdot, z)}\neq M_Y^D(\cdot, z)$.
\end{defn}

For any $z\in \partial_m D$, let $Y^{D, z}=(Y^{D, z}_t, \P^z_x)$
denote the $M^D_Y(\cdot, z)$-process, Doob's $h$-transform
of $Y^D$ with $h(\cdot)=M^D_Y(\cdot, z)$. The lifetime of $Y^{D, z}$
will be denoted by $\zeta$. It is known (see \cite{KW}) that
$\lim_{t\uparrow\zeta}Y^{D, z}_t=z$,  $\P^z_x$-a.s.
For $E\subset D$, let $T_E:= \inf\{t>0: Y^{D, z}_t\in E\}$.
It is proved in \cite[Satz 2.6]{Fol} that a set $E\subset D$
is minimally thin at $z\in \partial_m D$ if and only if
there exists $x\in D$ such that $\P^z_x(T_E<\zeta)\neq 1$.

We assume now that $z$ is a fixed point in $\partial D$  and the base point
$x_0$ of the Martin kernel $M^D_Y$
(cf.~\eqref{e:mk}) satisfies $x_0\in D\cap B(z, R)$ and $R/4<\delta(x_0)<R$.

The following criterion for minimal thinness has been proved
for a large class of symmetric L\'evy processes in \cite[Proposition 6.4]{KSV11}.
The proof is quite general and it works whenever (1) the cone of
excessive functions of the underlying process forms a balayage space,
and (2) the inequality in Corollary \ref{c:U-M} relating the
Green function and the Martin kernel of the processes is valid.
In particular, the proof works in the present setting.
For $E\subset D$, define
$$
E_n=E\cap \{x\in D:\, 2^{-n-1}\le |x-z|<2^{-n}\}\, ,\quad n\ge 1\, .
$$

\begin{prop}\label{p:minthin-criterion-1}
A set $E \subset D$ is minimally thin in $D$  at $z$ with respect
to $Y^D$ if and only if
$\sum_{n=1}^{\infty}R^{E_n}_{M^D_Y(\cdot, z)}(x_0)<\infty$.
\end{prop}

Let us fix $z\in \partial D$. Define $v(x)=U^D(x,x_0) \wedge 1$ so
that $v\in \SS^c(D)$. By Theorems \ref{t:gfe} and \ref{t:mke}
we see that for $x$ close to $z$,
$$
\frac{M^D_Y(x,z)}{v(x)}\asymp \frac{\phi'(|x-z|^{-2})}{|x-z|^{d+4}\phi(|x-z|^{-2})^2}
$$
with a constant depending on $z$ and $x_0$, but not on $x$.
By using Lemma \ref{l:properties-of-bf}(b),
we see that there exists a constant $c_1>0$ such that for large $n$,
$$
c_1^{-1} \frac{2^{n(d+4)}\phi'(2^{2n})}{\phi(2^{2n})^2}
\, v(x)\le M^D_Y(x,z)
\le c_1 \frac{2^{(n+1)(d+4)}\phi'(2^{2(n+1)})}{\phi(2^{2(n+1)})^2}\,
v(x)\, ,\quad x\in E_n\, .
$$
This implies that
$$
c_1^{-1} \frac{2^{n(d+4)}\phi'(2^{2n})}{\phi(2^{2n})^2}\, R^{E_n}_v
\le R^{E_n}_{M^D_Y(\cdot,z)}\le c_1 \frac{2^{(n+1)(d+4)}\phi'(2^{2(n+1)})}{\phi(2^{2(n+1)})^2}\, R^{E_n}_v\, .
$$
In particular,
\begin{align}\label{e:equivalence-1}
\sum_{n=1}^{\infty}R^{E_n}_{M^D_Y(\cdot, z)}(x_0)<\infty
\qquad \text{if and only if} \qquad \sum_{n=1}
\frac{2^{n(d+4)}\phi'(2^{2n})}{\phi(2^{2n})^2}\, R^{E_n}_v(x_0)<\infty\, .
\end{align}
Since $\wh{R}^{E_n}_{v}$ is a potential, there
is a measure $\lambda_n$ (supported by $\overline{E}_n$) charging no polar sets
such that
$\wh{R}^{E_n}_{v}=U^D \lambda_n$. Also, $\wh{R}^{E_n}_{v}=
v=U^D(\cdot, x_0)$ on $\overline{E}_n$ (except for a polar set,
and at least for large $n$), hence
\begin{eqnarray*}
\wh{R}^{E_n}_{v}(x_0)&=& U^D \lambda_n(x_0)=
\int_{\overline{E}_n}U^D(x_0,y)\, \lambda_n(dy)
=\int_{\overline{E}_n} v(y)\, \lambda_n(dy)\\
&=&\int_{\overline{E}_n}\wh{R}^{E_n}_{v}(y)\,
\lambda_n(dy)=\int_D \int_D U^D(x,y)\,
\lambda_n(dy)\, \lambda_n(dx)= \gamma_{v}(E_n)\, .
\end{eqnarray*}
We conclude from \eqref{e:equivalence-1} that
\begin{align}\label{e:equivalence-2}
\sum_{n=1}^{\infty}R^{E_n}_{M^D_Y(\cdot, z)}(x_0)<\infty
\qquad \text{if and only if} \qquad \sum_{n=1}
\frac{2^{n(d+4)}\phi'(2^{2n})}{\phi(2^{2n})^2}\, \gamma_v(E_n)<\infty\, .
\end{align}

Thus we have proved the following Wiener-type criterion for minimal thinness.

\begin{corollary}\label{c:minthin-criterion-1}
$E \subset D$ is minimally thin in $D$  at $z$ with respect to $Y^D$ if and only if
$$
\sum_{n=1}^{\infty} \frac{2^{n(d+4)}\phi'(2^{2n})}{\phi(2^{2n})^2}\, \gamma_{v}(E_n)<\infty.
$$
\end{corollary}
Now we state a version of Aikawa's criterion for minimal thinness.
\begin{prop}\label{p:aikawa-thinness}
Let $z\in \partial D$ and $E\subset D$,
let $\{Q_j\}$ be a Whitney decomposition of $D$ and let $x_j$ denote
the center of $Q_j$. The following are equivalent:
\begin{itemize}
\item[(a)] $E$ is minimally thin at $z$ with respect to $Y^D$;
\item[(b)]
\begin{align*}
\sum_{j: Q_j\cap B(z,1)\neq \emptyset}
\frac{v^2(x_j)\phi'(\mathrm{dist}(z,Q_j)^{-2})}{\mathrm{dist}(z,Q_j)^{d+4}\phi(\mathrm{dist}(z,Q_j)^{-2})^2}
\, \cp_D(E\cap Q_j)<\infty\, ;
\end{align*}
\item[(c)]
\begin{equation}\label{e:aikawa-thinness-1}
\sum_{j: Q_j\cap B(z,1)\neq \emptyset}
\frac{\mathrm{dist}^2(Q_j, \partial D)\phi'(
\mathrm{dist}(z,Q_j)^{-2})}{\mathrm{dist}(z,Q_j)^{d+4}\phi(
\mathrm{dist}(z,Q_j)^{-2})^2}\,
\cp_D(E\cap Q_j)<\infty\, .
\end{equation}
\end{itemize}
\end{prop}
\pf
By using Proposition \ref{p:quasi-additivity}, the proof is analogous to the proofs of
\cite[Proposition 6.6 and Corollary 6.7]{KSV11},
cf.~also \cite[Proposition 4.4]{MV}, therefore we omit the proof. \qed

\noindent
{\bf Proof of Theorem \ref{t:main}:} Assume that $E$ is
minimally thin at $z\in \partial D$. By Proposition
\ref{p:aikawa-thinness}, the series \eqref{e:aikawa-thinness-1}
converges. By Proposition \ref{p:comparable-measure}, the measure
$$
\sigma(A):=\int_A \phi(\delta(x)^{-2})\, dx\, ,\quad A\subset D\, ,
$$
is comparable to the capacity $\cp_D$ with respect to the
Whitney decomposition $\{Q_j\}$. Therefore
$$
\cp_D(E\cap Q_j)\ge c_1 \sigma(E\cap Q_j)=c_1 \int_E
\ind_{Q_j}(x)\phi(\delta(x)^{-2})\, dx\, .
$$
For $x\in Q_j$ we have that $\mathrm{dist}^2(Q_j,\partial D)
\asymp\delta(x)$ and $\mathrm{dist}(z,Q_j)\asymp |x-z| $. Therefore,
\begin{eqnarray*}
\infty&>&\sum_{j: Q_j\cap B(z,1)\neq \emptyset}
\frac{\mathrm{dist}^2(Q_j, \partial D)\phi'(\mathrm{dist}(z,Q_j)^{-2})}{\mathrm{dist}(z,Q_j)^{d+4}\phi(
\mathrm{dist}(z,Q_j)^{-2})^2}\, \cp_D(E\cap Q_j)\\
&\ge &c_2 \sum_{j: Q_j\cap B(z,1)\neq \emptyset}\int_E \frac{\delta(x)^2\phi'(|x-z|^{-2})
}{|x-z|^{d+4}
\phi(|x-z|^{-2})^2}\, \ind_{Q_j}(x)\phi(\delta(x)^{-2})\, dx\\
&=& c_2 \int_{E\cap B(z,1)}
\frac{\delta(x)^2\phi(\delta(x)^{-2})\phi'(|x-z|^{-2})}{|x-z|^{d+4} \phi(|x-z|^{-2})^2}\, dx \, .
\end{eqnarray*}

Conversely, assume that $E$ is a union of a subfamily of Whitney cubes of $D$.
Then $E\cap Q_j$ is either empty or equal to $Q_j$. Since
$\mathrm{Cap}_D(Q_j)\asymp \sigma(Q_j)=\int_{Q_j}
\phi(\delta(x)^{-2})\, dx$, we can reverse the first
inequality in the display above to conclude that
\begin{eqnarray*}
\lefteqn{\sum_{j: Q_j\cap B(z,1)\neq \emptyset}
\frac{\mathrm{dist}^2(Q_j, \partial D)\phi'(
\mathrm{dist}(z,Q_j)^{-2})}{\mathrm{dist}(z,Q_j)^{d+4}\phi(
\mathrm{dist}(z,Q_j)^{-2})^2}\, \cp_D(E\cap Q_j)}\\
&\le & c_3 \int_{E\cap B(z,1)}
\frac{\delta(x)^2\phi(\delta(x)^{-2})\phi'(|x-z|^{-2})}{|x-z|^{d+4}
\phi(|x-z|^{-2})^2}\, dx \, .
\end{eqnarray*}
\qed

Theorem \ref{t:main} will be now applied to study
minimal thinness of a set below the graph of a Lipschitz function.
We start by recalling Burdzy's result, cf.~\cite{Bur, Gar}:
Let $f:\R^{d-1}\to [0,\infty)$ be a Lipschitz function.
The set $A=\{x=(\wt{x},x_d) \in \H:\, 0<x_d\le f(\wt{x})\}$ is
minimally thin in $\H$ with respect to Brownian motion at $z=0$ if and only if
\begin{align}\label{c:criterion2f}
\int_{\{|\wt{x}|<1\}}f(\wt{x})|\wt{x}|^{-d}\, d\wt{x} <\infty\, .
\end{align}
It is shown recently in \cite{KSV6} that the same criterion
for minimal thinness is true for the subordinate Brownian motions
studied there. By using Theorem \ref{t:main} one can follow the
proof of \cite[Theorem 4.4]{KSV6} and show the Burdzy-type
criterion for minimal thinness in Proposition \ref{p:dahlberg2}.
In the proof we will need the following simple observation: For any $T>0$, we have for $t\in (0,T]$,
\begin{align}\label{e:integral}
\int_0^t r^2 \phi(r^{-2})\, dr \asymp t^3 \phi(t^{-2}),
\end{align}
Indeed, since  $r^2\phi(r^{-2})\le t^2\phi(t^{-2})$ for all $0<r\le t$, it follows that $\int_0^t r^2 \phi(r^{-2})\, dr \le t^3 \phi(t^{-2})$. On the other hand, since $\phi$ is increasing, $\int_0^t r^2 \phi(r^{-2})\, dr \ge   \phi(t^{-2})\int_0^t r^2 \, dr=  \frac{t^3}{3} \phi(t^{-2})$.

\begin{prop}\label{p:dahlberg2}
Assume that $d\ge 3$ and that $f:\R^{d-1}\to [0,\infty)$ is a Lipschitz function.
Suppose $D=\{x=(\wt{x},x_d)\in \R^d:\, x_d> h(\wt{x})\}$ is the domain above the graph of a bounded $C^{1,1}$ function $h$. Then the set
$$
A:=\{x=(\wt{x},x_d)\in \R^d:\, h(\wt{x}) <x_d\le f(\wt{x})+ h(\wt{x})\}
$$
is minimally thin in $D$ at $0$ with respect to $Y^D$ if and only if
\begin{align}\label{c:criterion2}
\int_{\{|\wt{x}|<1\}}\frac{f(\wt{x})^3\phi(f(\wt{x})^{-2})\phi'(|\wt{x}|^{-2})}
{|\wt{x}|^{d+4}\phi(|\wt{x}|^{-2})^2}\, d\wt{x} <\infty\, .
\end{align}
\end{prop}

\pf Without loss of generality we may assume that $f(\wt{0})=0$.
We first note that by the Lipschitz continuity of $f$, it follows
that $|\wt{x}|\le |x|\le c_1|\wt{x}|$ for $x=(\wt{x},x_d)\in A$.
Hence by Fubini's theorem we have
\begin{eqnarray}\label{e:main2-1}
\lefteqn{\int_A \frac{x_d^2 \phi(x_d^{-2})\phi'(|x|^{-2})}{|x|^{d+4}
\phi(|x|^{-2})^2}\, dx = \int_{|\wt{x}|< 1}d\wt{x} \int \ind_A(\wt{x}, x_d) \frac{x_d^2 \phi(x_d^{-2})\phi'(|x|^{-2})}{|x|^{d+4}
\phi(|x|^{-2})^2}\, dx_d}\nonumber \\
 &\asymp&  \int_{|\wt{x}|<1} \frac{\phi'(\wt{x})}{|\wt{x}|^{d+4}
 \phi(|\wt{x}|^{-2})^2}\, d\wt{x}\int_0^{f(\wt{x})} x_d^2 \phi(x_d^{-2})\, dx_d\nonumber\\
& \asymp&\int_{|\wt{x}|<1}\frac{f(\wt{x})^3\phi(f(\wt{x})^{-2})\phi'(|\wt{x}|^{-2})}
{|\wt{x}|^{d+2}\phi(|\wt{x}|^{-2})^2}\, d\wt{x},
\end{eqnarray}
where the last asymptotic relation follows from \eqref{e:integral}
with $T=\sup_{|\wt{x}|\le 1}f(\wt{x})$.
It follows from Theorem \ref{t:main} that if $A$ is minimally
thin in $D$ at $0$, then \eqref{c:criterion2} holds true.

For the converse,
let $\{Q_j\}$ be a Whitney decomposition of $D$ and define
$E= \cup_{Q_j\cap A\neq \emptyset}Q_j  $;
clearly $A\subset E$.
Let $Q_j^*$ be the interior of the double of $Q_j$ and
note that $\{Q_j^*\}$ has bounded multiplicity, say $N$.
Moreover, if $Q_j\cap A\neq \emptyset$,
then by the Lipschitz continuity of $f$ we have $|Q_j^*\cap A|\asymp |Q_j|$.
Moreover,
for $x \in Q_j^*$ we have
$|x|  \asymp \text{dist}(0, Q_j).$
Therefore
\begin{align}\label{e:main2-2}
&\int_A \frac{x_d^2 \phi(x_d^{-2})\phi'(|x|^{-2})}{|x|^{d+4}
\phi(|x|^{-2})^2}\, dx \le
\int_{E }\frac{x_d^2 \phi(x_d^{-2})\phi'(|x|^{-2})}{|x|^{d+4}
\phi(|x|^{-2})^2}\, dx\nonumber\\
&=\sum_{Q_j\cap A\neq \emptyset} \int_{Q_j  }\frac{x_d^2 \phi(x_d^{-2})\phi'(|x|^{-2})}{|x|^{d+4}
\phi(|x|^{-2})^2}\, dx \nn\\
 &\le  c_2\sum_{Q_j\cap A\neq \emptyset}  |Q_j^*\cap A|
 \frac{\mathrm{dist}^2(Q_j^*,D)\phi(\mathrm{dist}^{-2}(Q_j^*,D))\phi'(\mathrm{dist}^{-2}(0,Q_j))}
 {\mathrm{dist}^{d+4}(0,Q_j)\phi(\mathrm{dist}^{-2}(0,Q_j))^2}
\nn \\
& \le c_3\sum_{Q_j\cap A\neq \emptyset} \int_{Q_j^*\cap A}
\frac{x_d^2 \phi(x_d^{-2})\phi'(|x|^{-2})}{|x|^{d+4}
\phi(|x|^{-2})^2}\, dx
 \le c_3N \int_A \frac{x_d^2 \phi(x_d^{-2})\phi'(|x|^{-2})}{|x|^{d+4}
\phi(|x|^{-2})^2}\, dx\, .
\end{align}
If \eqref{c:criterion2} holds, then \eqref{e:main2-1} and
\eqref{e:main2-2} imply that $$\int_{E }\frac{x_d^2 \phi(x_d^{-2})\phi'(|x|^{-2})}{|x|^{d+4}
\phi(|x|^{-2})^2}\, dx     <\infty.$$ Hence, by Theorem \ref{t:main}, $E$ is minimally thin,
and thus $A$ is also minimally thin.
\qed

\begin{remark}\label{r:dahlberg2-2}{\rm
In case $d\ge 2$ and a \emph{bounded}  $C^{1,1}$ domain,
we can get an analog of Proposition \ref{p:dahlberg2}.
Let $z \in \partial D$ and choose a coordinate system
$CS$ with its origin at $z$
such that
$$
B(z, R)\cap D=\{ y= (\wt y, \, y_d) \mbox{ in } CS:
|y|< R, y_d > h (\wt y) \},
$$
where $h$ is a $C^{1,1}$-function $h: \bR^{d-1}\to \bR$
 satisfying
$h (\wt 0)= 0$.
 Let $f:\R^{d-1} \to [0,\infty)$ be a Lipschitz function and
 $$
A:=\{x=(\wt{x},x_d) \in D:\, |x|< R,  h(\wt{x}) <x_d\le f(\wt{x})+ h(\wt{x})\}.
$$
Then the set
is minimally thin in $D$ at $z \in \partial D$ with respect to $Y^D$ if and only if \eqref{c:criterion2} holds true.
}
\end{remark}


\section{Examples}\label{s:stable}

In this section we assume  $D$ is either a bounded $C^{1,1}$ domain in $\R^d$ or a half-space.
We first compare criteria for minimal thinness for three processes in $D$
related to the isotropic $\alpha$-stable process. The first process
is the killed isotropic $\alpha$-stable process $X^D$, $0<\alpha<2$,
that is a killed subordinate Brownian motion $X_t=W_{S_t}$
where $(S_t)_{t\ge 0}$ is an $(\alpha/2)$-stable subordinator.
The corresponding Laplace exponent is the function
$\phi(\lambda)=\lambda^{\alpha/2}$. The second process
is the subordinate killed Brownian motion $Y^D_t=W^D_{S_t}$
with the same  $(\alpha/2)$-stable subordinator. The third
process is the censored $\alpha$-stable process $Z^D$.
The process $Z^D$ is a symmetric Markov process with
Dirichlet form given by
$$
\CC(v,v)=\int_{D}\int_{D} (v(y)-v(x))^2 j(y-x)\, dy\, dx\, ,
$$
where $j(x)$ is the density of the L\'evy measure of the isotropic
$\alpha$-stable process. The censored stable process was introduced
and studied in \cite{BBC}. When $\alpha\in (1,2)$, $Z^D$ is transient
and converges to the boundary of $D$ at its lifetime.

Hardy's inequality for the Dirichlet form of $Z^D$ was obtained in \cite{CS, D}.
Let $G^D_Z$ be the Green function of $Z^D$.
If $D$ is a bounded $C^{1,1}$ domain,
sharp two-sided estimates on $G^D_Z$  were obtained in \cite{CK02}.
If $D$ is a half-space, say the upper half-space, then it follows from
\cite{BBC} that the censored $\alpha$-stable process in $D$ satisfies the following
scaling property: for any $c>0$, if $(Z^D_t)_{t\ge 0}$ is a censored $\alpha$-stable
process in $D$ starting from $x\in D$, then $(cZ^D_{t/c^{\alpha}})_{t\ge 0}$ is a
censored $\alpha$-stable process in $D$ starting from $cx$. Thus the transition density
$p^D_Z(t, x, y)$ of $Z^D$ satisfies the following relation:
$$
p^D_Z(t, x, y)=t^{-d/\alpha}p^D_Z(1, t^{-1/\alpha}x, t^{-1/\alpha}y), \qquad t>0, x, y\in D.
$$
Now using the short-time heat kernel estimates in \cite{CKS1} we immediately arrive
at the following global estimates:
$$
p^D_Z(t, x, y)\asymp t^{-\frac{d}{\alpha}}\left(1\wedge\frac{t^{1/\alpha}}{|x-y|}\right)^{d+\alpha}
\left(1\wedge\frac{\delta_D(x)}{t^{1/\alpha}}\right)^{\alpha-1}\left(1\wedge\frac{\delta_D(x)}{t^{1/\alpha}}\right)^{\alpha-1},
\qquad \mbox{on } (0, \infty)\times D\times D.
$$
Using the above estimates, one can easily get sharp two-sided estimates on $G_D$ from which one can easily show that 
$$\lim_{x \ni D \to \infty} \frac{G^D_Z(x,y)}{G^D_Z(z,y)}=0.$$
Sharp two-sided estimates on $G^D_Z$
give sharp two-sided estimates on the Martin kernel of $Z^D$.
The arguments in \cite{CK02} imply that
the finite part of the minimal Martin boundary of $D$ with respect to $Z^D$ and the finite
part of the  Martin boundary of $D$ with respect to $Z^D$ both coincide with
the Euclidean boundary $\partial D$ of $D$.

Based on these results,
one can follow the proof of \cite[Proposition 4.4]{MV} (which is an
analog of Proposition \ref{p:aikawa-thinness})  line by line and see
that the same results also hold when $D$ is a half-space.
Therefore the following holds.

\begin{prop}\label{p:aikawa-thinness-censored}
Let $\alpha\in (1,2)$ and $d \ge 2$.
Let $D$ be either a bounded $C^{1,1}$ domain in $\R^d$ or a half-space,
$z\in \partial D$, $E\subset D$, and let $x_j$ denote the center of $Q_j$.
Let $x_0\in D$ be fixed, $\cp^D$ be the capacity with respect to $Z^D$
and $v(x)=G^D_Z(x, x_0)\wedge 1$.
The following are equivalent:

\begin{itemize}
\item[(a)] $E$ is minimally thin at $z$;
\item[(b)]
\begin{align}\label{e:aikawa-thinness-2-censored}
\sum_{j: Q_j\cap B(z,1)\neq \emptyset}
\mathrm{dist}(z,Q_j)^{-d-\alpha+2} v(x_j)^2 \cp^D(E\cap Q_j)<\infty\, ;
\end{align}
\item[(c)]
\begin{align}\label{e:aikawa-thinness-censored}
\sum_{j: Q_j\cap B(z,1)\neq \emptyset}
\frac{\mathrm{dist}(Q_j,\partial D)^{2(\alpha-1)}}{\mathrm{dist}(z,Q_j)^{d+\alpha-2}} \, \cp^{D}(E\cap Q_j)<\infty\, .
\end{align}
\end{itemize}
\end{prop}

It is shown in \cite{MV} that the measure
$\sigma(A):=\int_A \delta(x)^{-\alpha}dx$ is comparable to
$\cp^D$ with respect to the Whitney decomposition.
Further, it follows from \cite[Theorem 1.1]{CK02} that
$v(x_j)\asymp  \mathrm{dist}(Q_j, \partial D)\asymp
\delta(x)^{2(\alpha-1)}$ for all $x\in Q_j$. With
this in hand one can use the argument in the proof of
Theorem \ref{t:main} to prove the following criterion
for minimal thinness with respect to the censored $\alpha$-stable process.

\begin{thm}\label{t:main-censored}
Assume that $\alpha\in (1,2)$.
Let $D$ be either a bounded $C^{1,1}$ domain in $\R^d$ or a half-space, $d\ge 2$, and
let $E$ be a Borel subset of $D$.

\noindent
(1)
If $E$ is  minimally thin in $D$  at $z\in\partial D$ with respect to $Z^D$, then
$$
\int_{E\cap B(z, 1)}
\frac{\delta(x)^{\alpha-2}}{|x-z|^{d+\alpha-2}}\, dx < \infty\, .
$$

\noindent
(2)
Conversely, if $E$ is the union of a subfamily of Whitney cubes of $D$ and is
not minimally thin in $D$  at $z\in\partial D$ with respect to $Y^D$, then
$$
\int_{E\cap B(z, 1)}
\frac{\delta(x)^{\alpha-2}}{|x-z|^{d+\alpha-2}}\, dx = \infty\, .
$$
\end{thm}

Note that for $X^D$ the integral in the criterion for minimal thinness is
$$
\int_{E\cap B(z, 1)}
\frac{1}{|x-z|^d}\, dx\, ,
$$
while for $Y^D$ the corresponding integral becomes
$$
\int_{E\cap B(z, 1)}
\frac{\delta(x)^{2-\alpha}}{|x-z|^{d+2-\alpha}}\, dx\, .
$$

\begin{corollary}\label{c:comparison}
Let $D$ be either a bounded $C^{1,1}$ domain in $\R^d$ with $d \ge 2$ or a half-space with $d \ge 3$.
Let $E$ be the union of a subfamily of Whitney cubes of $D$
and $z\in \partial D$.

\noindent
(i) Let $1<\alpha <2$. If $E$ is minimally thin at $z$ with
respect to  $Z^D$, then it is minimally thin at $z$ with respect to $X^D$.

\noindent
(ii)  Let $0< \alpha <2$. If $E$ is minimally thin at $z$ with
respect to  $X^D$, then it is minimally thin at $z$ with respect to $Y^D$.

\noindent
(iii) Let $1<\alpha_1 \le \alpha_2<2$. If $E$ is minimally thin at
$z$ with respect to the $\alpha_1$-stable censored process, then
it is minimally thin at $z$ with respect to the $\alpha_2$-stable
censored process.

\noindent
(iv) Let $0<\alpha_1 \le \alpha_2<2$. If $E$ is minimally thin at $z$
with respect to $Y^D$ with index $\alpha_2$, then it is minimally thin
at $z$ with respect to $Y^D$ with index $\alpha_1$.
\end{corollary}
\pf  All statements follow easily from criteria in Theorems
\ref{t:main} and \ref{t:main-censored} together with the
observation that since $\delta(x)\le |x-z|$,
$$
\left(\frac{\delta(x)}{|x-z|}\right)^{2-\alpha}\le 1 \le
\left(\frac{\delta(x)}{|x-z|}\right)^{\alpha-2}\, .
$$
\qed

A criterion for minimal thinness of a set below the graph of
a Lipschitz function with respect to the censored stable
process is given in the following result which can be
proved in the same way as Proposition \ref{p:dahlberg2}.
\begin{prop}\label{p:dahlberg2-censored}
Let $\alpha\in (1,2)$. Assume that $f:\R^{d-1}\to [0,\infty)$
is a Lipschitz function. Suppose that $D=\{x=(\wt{x},x_d)\in \R^d:\, 0<x_d \}$.
Then the set
$$
A:=\{x=(\wt{x},x_d)\in D:\, 0<x_d\le f(\wt{x})\}
$$
is minimally thin in $D$ at $0$ with respect to $Z^D$ if and only if
\begin{align}\label{c:criterion2-Z}
\int_{\{|\wt{x}|<1\}}\frac{f(\wt{x})^{\alpha-1}}{|\wt{x}|^{d+\alpha-2}}d\wt{x} <\infty\, .
\end{align}
\end{prop}

In case of $X^D$, the criterion reads
\begin{align}\label{c:criterion2-X}
\int_{\{|\wt{x}|<1\}}\frac{f(\wt{x})}{|\wt{x}|^{d}}d\wt{x} <\infty\, ,
\end{align}
while for $Y^D$ with $d \ge 3$,  \eqref{c:criterion2} becomes
\begin{align}\label{c:criterion2-Y}
\int_{\{|\wt{x}|<1\}}\frac{f(\wt{x})^{3-\alpha}}{|\wt{x}|^{d+2-\alpha}}d\wt{x} <\infty\, .
\end{align}

\begin{example}\label{ex:log}{\rm
Let $d \ge 3$ and $D=\{x=(\wt{x},x_d)\in \R^d:\, 0<x_d \}$,
$f:\R^{d-1}\to [0,\infty)$ a Lipschitz function
and put $A:=\{x=(\wt{x},x_d)\in D:\, 0<x_d\le f(\wt{x})\}$.

\noindent
(1) If  $f(\wt{x})=|\wt{x}|^{\gamma}$ with $\gamma\ge 1$, then an
easy calculation shows that all three integrals in
\eqref{c:criterion2-Z}-\eqref{c:criterion2-Y} are finite
if and only if $\gamma>1$. Thus, for all three processes,
$A$ is minimally thin at $z=0$ if and only if $\gamma>1$.

\noindent
(2) Let $f(\wt{x})=|\wt{x}|\big(\log(1/|\wt{x}|)\big)^{-\beta}$,
$\beta\ge 0$. Then $f$ is Lipschitz. By use of
\eqref{c:criterion2-Z}-\eqref{c:criterion2-Y} it follows
easily that $A$ is minimally thin at $z=0$
\begin{align*}
& \text{with respect to }Z^D \text{ if and only if } \beta>\frac{1}{\alpha-1}\, ,\\
& \text{with respect to }X^D \text{ if and only if } \beta>1\, ,\\
& \text{with respect to }Y^D \text{ if and only if } \beta>\frac{1}{3-\alpha}\, .
\end{align*}
Since $1<1/(3-\alpha)$ for $\alpha\in (0,2)$ and $1<1/(\alpha-1)$
for $\alpha\in (1,2)$ this is in accordance with Corollary \ref{c:comparison}.
By choosing $\beta$ and $\alpha$ appropriately, we conclude that none of
the converse in Corollary \ref{c:comparison} holds true.
}
\end{example}

We conclude this paper with an example about minimal thinness with respect to
subordinate killed Brownian motion in the half-space via geometric stable
subordinators.
We define $L_1(\lambda)=\log\lambda$, and for $n\ge 2$ and $\lambda>0$ large enough,
$L_n(\lambda)=L_1(L_{n-1}(\lambda))$.
Applying Proposition \ref{p:dahlberg2}, we can easily check the following.

\begin{example}\label{ex:log2}{\rm
Let $d \ge 3$ and
$\alpha\in (0, 1]$.
Suppose that $D=\{x=(\wt{x},x_d)\in \R^d:\, 0<x_d \}$ and $Y^D$ is the
subordinate killed Brownian motion in $D$ via a subordinator with Laplace
exponent $\log(1+\lambda^{\alpha})$.
Assume that $f:\R^{d-1}\to [0,\infty)$ a Lipschitz function
and define $A:=\{x=(\wt{x},x_d)\in D:\, 0<x_d\le f(\wt{x})\}$.

\noindent
(1) Let $f(\wt{x})=|\wt{x}|\big(
L_1(1/|\wt{x}|)\big)^{-\beta}$ with
$\beta\ge 0$. Then $A$ is minimally thin at $z=0$ with respect to $Y^D$  if and only if  $\beta> 0$.

\noindent
(2)  Let $n\ge 2$ and $f(\wt{x})=|\wt{x}|(L_2(1/|\wt{x}|)\cdots L_n(1/|\wt{x}|)\big)^{-1/3}
\big(L_{n+1}(1/|\wt{x}|)\big)^{-\beta}$ with
$\beta\ge 0$. Then $A$ is minimally thin at $z=0$ with respect to $Y^D$  if and only if  $\beta> 1/3$.
}
\end{example}

\bigskip
\noindent
{\bf Acknowledgements:} We are grateful to the referee for the insightful comments
on the first version of this paper.

\vspace{.1in}
\begin{singlespace}


\end{singlespace}

\end{doublespace}
\vskip 0.1truein

\parindent=0em

{\bf Panki Kim}

Department of Mathematical Sciences and Research Institute of Mathematics,

Seoul National University, Building 27, 1 Gwanak-ro, Gwanak-gu Seoul 151-747, Republic of Korea

E-mail: \texttt{pkim@snu.ac.kr}

\bigskip

{\bf Renming Song}

Department of Mathematics, University of Illinois, Urbana, IL 61801,
USA

E-mail: \texttt{rsong@math.uiuc.edu}

\bigskip

{\bf Zoran Vondra\v{c}ek}

Department of Mathematics, University of Zagreb, Zagreb, Croatia, and \\
Department of Mathematics, University of Illinois, Urbana, IL 61801,
USA

Email: \texttt{vondra@math.hr}
\end{document}